\documentclass[10pt,a4paper,oneside]{article}

\usepackage[margin=3cm,top=3cm,bottom=3cm]{geometry}
\usepackage[british,UKenglish]{babel}
\usepackage{csquotes}
\usepackage[final]{pdfpages}

\usepackage{amsthm}
\usepackage{amsfonts}
\usepackage{amsmath}
\usepackage{amssymb}
\usepackage{mathtools}
\usepackage{centernot}
\usepackage{stmaryrd}

\usepackage{pythonhighlight}

\usepackage{tikz}
\usepackage{tikz-cd}
\usepackage{genyoungtabtikz}
\usetikzlibrary{decorations.pathreplacing,angles,quotes,calc,positioning}

\usepackage{marginnote}
\usepackage{soul}
\usepackage[normalem]{ulem}
\mathtoolsset{showonlyrefs}
\usepackage{comment}
\usepackage{enumitem}

\usepackage{float}
\usepackage{import}
\usepackage{xifthen}
\usepackage{pdfpages}
\usepackage{transparent}
\usepackage{afterpage}
\usepackage{adjustbox}

\graphicspath{ {./figures/} }

\usepackage{hyperref}
\hypersetup{colorlinks=true,citecolor=blue,urlcolor =black,linkbordercolor={1 0 0}}
\usepackage[style=alphabetic,maxbibnames=99]{biblatex}
\newtheorem{thm}{Theorem}[section]
\newtheorem{conj}{Conjecture}
\newtheorem*{thm*}{Theorem}

\newtheorem{lemma}[thm]{Lemma}
\newtheorem{prop}[thm]{Proposition}

\newtheorem{prop-def}[thm]{Proposition/Definition}

\theoremstyle{definition}
\newtheorem{definition}[thm]{Definition}
\newtheorem{ex}[thm]{Example}

\theoremstyle{remark}
\newtheorem{rem}[thm]{Remark}

\newcommand{\eeuler}{\mathsf{e}}

\newcommand{\V}{\mathsf{V}}
\newcommand{\vir}{\mathrm{vir}}

\newcommand{\C}{\mathbb{C}}
\newcommand{\E}{\mathbb{E}}
\newcommand{\Z}{\mathbb{Z}}

\newcommand{\T}{\mathbf{T}}

\newcommand{\N}{\mathbb N}
\newcommand{\PS}{\mathbb P}

\newcommand{\coker}{\operatorname{coker}}

\newcommand{\length}{\operatorname{length}}
\newcommand{\fix}{\operatorname{fix}}
\newcommand*\spvertund[1]{\vrule width0pt height0pt depth#1\relax}

\DeclareFontFamily{U}{mathc}{}
\DeclareFontShape{U}{mathc}{m}{it}%
{<->s*[1.03] mathc10}{}
\DeclareMathAlphabet{\mathcal}{U}{mathc}{m}{it}

\DeclareMathOperator{\sHom}{\mathcal{H\mkern-3mu om}}
\DeclareMathOperator{\sExt}{\mathcal{E\mkern-3mu xt}}

\title{On the Virtual Euler Characteristic of the \\Moduli Space of Stable Pairs on Surfaces}

\author{Ana Pavlakovi\'c}

\begin{document}
	\maketitle

\begin{abstract}
	We study the stable pair theory on toric surfaces and determine the virtual tangent space over the fixed point loci. Further, we present a program to compute the virtual Euler characteristic, illustrated by the case of the projective plane. 
    As an application, conjectures regarding rationality and symmetry 
    are supported by verification of a special case. 
\end{abstract}

	\pagenumbering{arabic}
	
	\section{Introduction}

Theories of virtual invariants and their relations to each other have played an important role in the study of enumerative geometry \cite{pt_13/2}.  Three are particularly important:  the theory of stable pairs, developed by Pandharipande and Thomas \cite{pt_derived} (PT-theory), Donaldson-Thomas theory \cite{thomas} (DT-theory), and Gromov-Witten theory (GW-theory).  In the case of toric $3$-folds, these are all equivalent \cites{pp,moop,mnop1}. 
 
 So far, PT-theory has mainly been applied to the case of threefolds.  However, from PT-theory we have learned valuable information about surfaces.  For instance, in \cite{pt_nested} an  equivalence  was shown between the moduli space of stable pairs on surfaces and a specific nested Hilbert scheme, which  we use to apply PT-theory to the study of surfaces. Moreover, important methods were developed, e.g., box configurations were introduced in \cite{mnop1}, where it was conjectured that DT-theory and GW-theory are equivalent.  We use box configurations in the proof of our main result, which is the computation of the virtual tangent space of the moduli space of stable pairs on a toric surface. 

\begin{thm}\label{prop:main conjecture is equiv to loc}
	For a non-singular toric surface $S$ and  a homology class $\beta\in H_2(S;\Z)$, 
	\begin{equation}
		\tilde{Z}^{S}_{\beta} (q\mid s_1,s_2)  = \sum_{m\geq 0} \sum_{[\vec{d}] = \beta}   \sum_{\substack{\vec{\lambda}\subset \vec{d} \\|\vec{\lambda}| = m}} q^{m}  \prod_{\alpha\beta} \frac{c_\bullet}{\eeuler}\left(\mathsf{E}_{\alpha\beta}\right) \prod_{\alpha} \frac{c_\bullet}{\eeuler} \left(\mathsf{F}_\alpha + \mathsf{G}_\alpha\right),
	\end{equation}
	where 
	\begin{equation}
		T^\vir = \sum_{\alpha\beta} \mathsf{E}_{\alpha\beta} + \sum_\alpha (\mathsf{F}_\alpha + \mathsf{G}_\alpha)
	\end{equation}
	and $\mathsf{G}_\alpha, \mathsf{F}_\alpha$, and $\mathsf{E}_{\alpha\beta}$ are defined as in \eqref{eq:Galpha}, \eqref{eq:vertex contribution}, and \eqref{eq:edge contribution}, respectively.
\end{thm}

Theorem \ref{prop:main conjecture is equiv to loc} is proven in Section \ref{sec:virtual tangent space}.


 A first application of Theorem \ref{prop:main conjecture is equiv to loc} is regarding a rationality conjecture, which we next motivate. 
 Let $S$ be a smooth projective surface and let $S^{[n]}$ be its Hilbert scheme of $n$ points. G\"ottsche computed the generating series of the topological Euler characteristics as the function
	\begin{equation}\label{eq:goettsche}
		\sum_{n\geq 0} q^n \eeuler(S^{[n]}) = \left(q^{-\frac{1}{24}}\eta(q)\right)^{-\eeuler(S)},
	\end{equation}
	where $\eta$ denotes the Dedekind eta function  \cite{goettsche}. The rationality of the series makes computing the virtual Euler characteristics straightforward.
The motivation for studying the rationality of the series for stable pairs on surfaces is twofold.

First, Oprea and Pandharipande introduced and studied the generating series of the virtual Euler characteristic of the Quot scheme of a surface in \cite{op}:
	\begin{equation}\label{eq:Z QUOT}
		Z_{S,N,\beta} \coloneq \sum_{n\in\Z}q^n \eeuler^\vir (\operatorname{Quot}(\C^N, \beta,n)).
	\end{equation}
	The virtual Euler characteristic of a space $X$ with a $2$-term perfect obstruction theory,
	\begin{equation}
		\eeuler^\vir (X) = \int_{[X]^\vir} c_\bullet(T^\vir X),
	\end{equation}
	defined by \cite{fg}, arises from the virtual tangent space of the canonical obstruction theory of the space. 
	In \cite{op}, Oprea and Pandharipande conjectured that \eqref{eq:Z QUOT} is the Laurent expansion of a rational function. In \cite{jop}, Johnson, Oprea, and Pandharipande showed that if $S$ is simply connected and $N=1$, in which case 
	\begin{equation}\label{eq:quot as product}
		\operatorname{Quot} \simeq S^{[n]} \times \operatorname{Hilb}_\beta,
	\end{equation}
	the series $Z_{S,1,\beta}$ is the expansion of a rational function. In \cite{lim}, Lim showed the rationality conjecture for surfaces with $p_g >0$ for any $N$.

Second, in 2013, Pandharipande and Pixton showed that  the capped descendent series of stable pairs on threefolds is a rational function \cite{pp}. 
The moduli space of stable pairs was first introduced in \cite{pt_derived} as a means to study the curve counting theory of three-folds  \cite{pt_vertex}. 
	The moduli space can be thought of as parametrising curves and points contained in the curves. In the surface case, it is the natural extension of the Hilbert scheme of points and the Hilbert scheme of divisors. On one hand, one can study the product of the two spaces \eqref{eq:quot as product}, as done in \cite{op} and \cite{jop}, but one can also ask whether the series is rational if one restricts the points to be contained in the curves. More precisely, we study the rationality of the following series:
    \begin{equation}
        Z_{\beta}^S (q) = \sum_{n \in \N} q^n \eeuler^\vir ( P_n(S,\beta)),
    \end{equation}
    where $S$ is a projective surface, $\beta \in H_2(S)$ is a homology class, and $P_n(S,\beta)$ is the moduli space of stable pairs on $S$ with fixed Chern class $\beta$ and Euler characteristic $n$.
    For $S = \PS^2$ and $\beta = 1$, we explicitly determine the rational function in Section \ref{sec:d=1}:
	\begin{thm}\label{thm:P2,d=1}
		The series 
    \begin{equation}
        Z^{\PS^2}_{d=1} (q) = \sum_{n\in \N} q^n \eeuler^\vir(P_n(\PS^2, 1))
    \end{equation}
    is the Laurent expansion of 
    \begin{equation}
        \frac{3}{(1-q)^2}.
    \end{equation}
	\end{thm}

 Due to further explicit computations for $\mathbb{P}^2$ in low degree, which lead to concrete conjectures for rational functions (see Section \ref{sec:examples}), we believe that this result generalises to the following conjecture:
	\begin{conj}\label{conj:main conjecture}
		For a non-singular projective surface $S$ and a homology class $\beta\in H_2(S;\Z)$,  
		\begin{equation}
			Z^{S}_{\beta} (q) = \sum_{n\in\Z} q^n \int_{[P_n(S,\beta)]^\vir} c_\bullet \left(T^\vir_{P_n(S,\beta)}\right)
		\end{equation}
		is the Laurent expansion of a rational function of $q$.\footnote{We also believe the more general series with descendent insertions to be rational.}
		Moreover, the series viewed as a rational function satisfies the symmetry 
		\begin{equation}
			Z_{\beta}^{S}(q) = Z_{\beta}^{S}\left(\frac{1}{q}\right)
		\end{equation}
		and has possible poles only at $0$ and $1$. 
	\end{conj}

Studying toric surfaces turns out to be a good approach, because of their combinatorial nature.
A key tool for studying moduli spaces over toric varieties is virtual localisation (see Section \ref{sec:fixed points} for a detailed introduction). 
In the case of toric varieties, Theorem \ref{thm:P2,d=1} generalises to the following equivariant statement:
	\begin{conj}
		The series
		\begin{equation}\label{eq:Z tilde}
			\tilde{Z}^{S}_{\beta} (q\mid s_1,s_2) \coloneq \sum_{m\geq 0} \sum_{[\vec{d}] =\beta}   \sum_{\substack{\vec{\lambda}\subset \vec{d} \\|\vec{\lambda}| = m}} q^{m} \frac{c_\bullet}{\eeuler}  \left(T^\vir\right) 
		\end{equation}
		is the Laurent expansion of a rational function in $q,s_1,s_2$. 
	\end{conj}
By \eqref{eq:equivariant and non-equivariant series relation}, it implies Conjecture \ref{conj:main conjecture} for compact toric surfaces.
	
	
\subsection{Non-Rationality of Local Partition Function}

Studying the global partition function proves to be very difficult, so a natural approach is to first study the local partition functions. That is, 
\begin{equation}
		\tilde{Z}^{\C^2}_{(d_1,d_2)} (q\mid s_{1},s_{2}) \coloneqq \sum_{\lambda\subset (d_{1}, d_{2})} q^{|\lambda|} \frac{c_\bullet}{\eeuler} \left( \mathsf{G}_{\lambda} \right)
	\end{equation}
for
\begin{align}
		\tilde{Z}_\beta^S(q\mid s_1,s_2) &= \sum_{[\vec{d}] =\beta} \prod_{\alpha\beta} \frac{c_\bullet}{\eeuler} (\mathsf{E}_{\alpha\beta}) \prod_\alpha   \frac{c_\bullet}{\eeuler} (\mathsf{F}_\alpha) \sum_{\vec{\lambda}\subset\vec{d}} q^{|\vec{\lambda}|} \prod_\alpha \frac{c_\bullet}{\eeuler} (\mathsf{G}_\alpha)\\
		&=  \sum_{[\vec{d}] =\beta} \prod_{\alpha\beta} \frac{c_\bullet}{\eeuler} (\mathsf{E}_{\alpha\beta}) \prod_{\alpha} \frac{c_\bullet}{\eeuler}(\mathsf{F}_\alpha) \prod_\alpha Z_{(d_{\alpha\beta_1},d_{\alpha\beta_2})}^{U_\alpha} (q\mid s_{\alpha\beta_1},s_{\alpha\beta_2}).
\end{align}
Concretely, in the case where $S = \PS^2$ and $d=1$, 
\begin{multline}
		\sum_{\vec{\lambda}\subset (0,1,1)} q^{|\vec{\lambda}|}\frac{c_\bullet}{\eeuler} (\mathsf{G}_\alpha) \frac{c_\bullet}{\eeuler} (\mathsf{G}_\beta)  \frac{c_\bullet}{\eeuler} (\mathsf{G}_\gamma) \\
		= \tilde{Z}_{(0,1)}^{\C^2} (q\mid -s_1,s_2-s_1)\cdot \tilde{Z}_{(1,1)}^{\C^2} (q\mid s_1-s_2,-s_2) \cdot \tilde{Z}_{(1,0)}^{\C^2} (q\mid s_2,s_1).
	\end{multline}
However, it turns out that these are not expansions of rational functions. It fails even in the simplest case \eqref{eq:(0,1)}, 
\begin{equation}
		\tilde{Z}_{(0,1)}^{\C^2} (q\mid s_{1},s_{2}) = \sum_{b\geq 0}q^b \prod_{i=1}^b \frac{1-is_{1}}{-is_{1}}  = (1-q)^{\frac{1}{s_1} - 1} \notin \mathbb{Q}(s_{1})(q).
\end{equation}

Despite the non-rationality of the local partition function, we believe that the non-rationality of controlled, that is, that the non-rational factors cancel out after multiplying the local partition functions to form the global partition function.
	
\subsection{Further Directions}

Capped localisation introduced in \cite{moop} was used in \cite{pp} to show rationality of the global partition function for toric 3-folds. We presume that using capped localisation one can show rationality of the global partition function for surfaces, too. Here, in the surface case, one need only consider the capped descendent $1$- and $2$-leg vertex. Furthermore, in \cite{pp} they prove the rationality of the descendent partition function for stable pairs on nonsingular toric 3-folds. We believe that a generalisation of the conjecture to include descendents also holds.

\subsection{Plan of the Paper}
	
	In Section \ref{sec:stable pairs} we recall the notion of stable pairs on surfaces and its deformation theory. Moreover, we give a geometric description via the nested Hilbert scheme. 
	In Section \ref{sec:fixed points} we determine the fixed points of the induced torus action on the moduli space of stable pairs. To this end, we introduce Young diagrams and partitions and relate them to the fixed points. 
	In Section \ref{sec:virtual tangent space}, we use this description to explicitly compute the virtual tangent space of $P_{n}(S,d)$ and then establish that it has no fixed part. This is needed for the class $\frac{c_\bullet}{\eeuler} (T^\vir)$ to be well-defined. Here we use the technique of regrouping used in \cite{mnop1} and \cite{pt_vertex} to describe the virtual tangent space purely in terms of Laurent polynomials. 
	We can then, in Section \ref{sec:chern euler}, compute $\frac{c_\bullet}{\eeuler} (T^\vir)$. 
	Finally, in Section \ref{sec:examples}, we provide explicit examples of the virtual Euler characteristic for low $d$ both locally and globally.

	\subsection{Notation}
	We fix some notation.
	\begin{enumerate}[label =(\roman*)]
		\item The alternating $\operatorname{Ext}$-group for coherent sheaves $F,G$ over a scheme $X$
		\begin{equation}
			\sum_{i=0}^\infty (-1)^i \operatorname{Ext}^i_X(F,G)
		\end{equation}
		is denoted by  $\chi(F,G)$.
		\item Note that as we are working with toric localisation and are computing classes in $K^\T(\mathrm{pt})$, where $\T$ denotes the torus with coordinates $t_1,t_2$, we have to take the trace of elements in $K$-theory to be able to write them in terms of the torus coordinates $t_1,t_2$. The trace is the isomorphism between $K^\T(\mathrm{pt})$ and the algebra of Laurent polynomials in $t_1$ and $t_2$ and shall be omitted from the notation.
		\item For a torus invariant subscheme $X$ of $\C^2_{x_1,x_2}$ with ideal sheaf $\mathcal{I}_X$, we denote 
		\begin{equation}
			Q_{X}(t_{1},t_{2}) \coloneqq \sum_{\substack{(i,j)\in\Z^2_{\geq0},\\x_1^ix_2^j \notin \mathcal{I}_X}} t_{1}^{i} t_{2}^{j}.
		\end{equation}
		If $t_{1},t_{2}$ are clear from the context, we omit them and simply write $Q_X$.
		\item Let $\overline{Q}_X(t_{1},t_{2}) \coloneqq Q_X(t_{1}^{-1},t_{2}^{-1})$.  For the alternating $\operatorname{Ext}$-group we use $\chi(F)^*$ instead of $\overline{\chi(F)}$. 
	\end{enumerate}

\subsection{Acknowledgements}

The author would like to thank M. Moreira for proposing the problem and, jointly with W. Lim, advising the author's Master's thesis from which this paper stems. In particular, \eqref{eq:(1,1)} is due to M. Moreira and W. Lim contributed to \eqref{eq:S=P2, d=1}. Furthermore, discussions with  A. Iribar Lopez played an important role. A. Gibney provided valuable input during the writing of this manuscript. Lastly, the author would like to thank R. Pandharipande, without whom this project would not have been possible.





	\section{Stable Pairs on Surfaces}\label{sec:stable pairs}
	
	\subsection{Stable Pairs}
	
	\begin{definition}
		A \emph{stable pair} $(F,s)$ on a surface $S$ is a coherent sheaf $F$ on $S$ together with a section $s \colon\mathcal{O}_S \to F$ such that 	
		\begin{enumerate}[label=(\roman*)]
			\item $F$ is pure of dimension $1$,
			\item $\dim \coker s = 0$.
		\end{enumerate}
	\end{definition}
	
	\begin{rem}
		To every such stable pair we can associate the exact sequence
		\begin{equation}
			0\to \ker \to \mathcal{O}_S \xrightarrow{s} F \to \coker \to 0.
		\end{equation}
	\end{rem}
	and  the following numerical data:
	\begin{enumerate}[label =  (\roman*)]
		\item $c_1(F) = \beta \in H^2(S;\Z)$,
		\item $\chi(S,F) = n \in\Z$.
	\end{enumerate}
	If we fix a cohomology class $\beta\in H^2(S;\Z)$ and an integer $n\in\Z$, we can consider the moduli space of all stable pairs such that $c_1(F) = \beta$ and $\chi(S,F) = n$.  
	This moduli space is denoted by $P_n(S,\beta)$ and is constructed in \cite{potier}.

	\subsection{Nested Hilbert Schemes}\label{sec:nested hilbert schemes}
	Another way to view $P_n(S,\beta)$, is to regard it as a nested Hilbert scheme: 
	
	\begin{definition}
		The moduli space $S^{[m_1,m_2]}_\beta$, called a \emph{nested Hilbert scheme} parametrises tuples $(Z_1,Z_2,D)$ such that $Z_1,Z_2$ are $0$-dimensional subschemes, $D$ is a pure $1$-dimensional subscheme (in this case a divisor), and $I_{Z_1}(-D) \subset I_{Z_{2}}$. We require $\length(Z_1) = m_1$, $\length(Z_2) = m_2$, and $[D] = \beta$.
	\end{definition}	
	We will only concern ourselves with the case $m_1=0$, as this is equivalent to the moduli space of stable pairs. 	
	\begin{prop}[{\cite[Proposition B.8]{pt_nested}}]\label{prop:moduli correspondence}
		The relative Hilbert scheme  $S_\beta^{[0,m]}$  is isomorphic to the moduli space of stable pairs $P_{n}(S,\beta)$, 
		where 
		\begin{equation}
			m = n + \frac{\beta(\beta+K_S)}{2},
		\end{equation}
		via the mapping 
		\begin{equation}
			(Z,D) \mapsto (\mathcal{O}_S \to \mathcal{O}_D(Z)).
		\end{equation}		
	\end{prop}	
	
	\begin{rem}\label{rem:m and n}
		Other examples are:
		\begin{enumerate}[label = (\arabic*)]
			\item If  $\beta\neq0$, $m_1\neq 0$, and $m_2= 0$, then $S^{[m_1,0]}_\beta = \operatorname{Quot}(n,\beta)$ is the Quot scheme of $S$ (where $m_1$ and $n$ are related as above). 
			\item If  $\beta \neq0$ and $m_1, m_2 = 0$, then $S^{[0,0]}_\beta = \operatorname{Hilb}_\beta (S)$ is the Hilbert scheme of divisors of class $\beta\in H^2(S;\Z)$.
			\item If $\beta = 0$ and $m_1 = m_2 = m$, then $S^{[m,m]}_{\beta=0} = S^{[m]}$ is the Hilbert scheme of $m$ points. 
		\end{enumerate}
	\end{rem}

	\subsection{Universal Pair and Deformation Theory}
	
	The argument for the existence of the universal pair is the same as in the threefold case treated in \cite[Section 2.3]{pt_derived}. 
	The moduli space is the GIT quotient of a subset of the product of a Quot scheme and a Grassmannian. There is a universal sheaf pulled back from the Quot scheme, with a universal section over the product, see \cite{potier}. 
	Over the stable locus there are no stabilisers, so by Kempf's lemma the universal pair descends to $P_n(S,\beta)$. 
	We denote it by 
	\begin{equation}
		\mathcal{O}_{S\times P_n(S,\beta)} \to\mathcal{F}
	\end{equation}
	and by
	\begin{equation}
		\mathsf{pr}_S \colon S\times P_n(S,\beta) \to S,
	\end{equation}
	\begin{equation}
		\mathsf{pr}_{P} \colon S\times P_n(S,\beta) \to P_n(S,\beta)
	\end{equation}
	we denote the projections onto the first and second factors. 
	Under the isomorphism $P_n(S,\beta)\cong S_\beta^{[0,m]}$, we have that 
	\begin{equation}
		\mathcal{F} \simeq \mathcal{O}_{\mathcal{D}_\beta}(\mathcal{Z}),
	\end{equation}
	where  $\mathcal{Z}$ is the universal $0$-dimensional subscheme  and  $\mathcal{D}_\beta$ the universal divisor.
	
	Let us first consider the moduli space as parametrising stable pairs. We then have a perfect obstruction theory by \cite[Theorem 2]{lin}: 
	\begin{equation}\label{eq:pot stable pairs} 
		\E _\bullet = R\sHom_{\mathsf{pr}_P}([\mathcal{O}_{S\times P_n(S,\beta)}\xrightarrow{s}\mathcal{F}],\mathcal{F}) \in D_{\mathrm{perf}}^{[0,1]}(P_n(S,\beta)).
	\end{equation} 
	
	On the other hand, by \cite[Proposition 2.5]{gsy} we also have: 
	\begin{equation}\label{eq:pot nested Hilbert scheme}
		\E_\bullet = \operatorname{Cone} \Big( R\sHom_{\mathsf{pr}_P}(\mathcal{I}_{\mathcal{Z}},\mathcal{I}_\mathcal{Z})\to R\sHom_{\mathsf{pr}_P}(\mathcal{O},\mathcal{I}_\mathcal{Z}(\mathcal{D}_\beta))\Big)
	\end{equation}
	As the moduli space has a $2$-term perfect obstruction theory, it has a virtual fundamental class.
	
	\section{$\T$-Fixed Points}\label{sec:fixed points}
	
	We study the virtual Euler characteristic series via localisation.

	\begin{thm}[{\cite{gp}}]\label{thm:virtual localisation}
		Let $X$ be a scheme with a $\C^\times$-action and a $\C^\times$-equivariant perfect obstruction theory. Then 
		\begin{equation}
			[X]^\vir = \iota_* \sum \frac{[X_i]^\vir}{e(N_i^\vir)} \in A_*^{\C^\times}(X)\otimes\mathbb{Q}\left[s,\frac{1}{s}\right].
		\end{equation}
		Here the $X_i$ denote the connected components of the fixed loci of $X$ under the $\C^\times$-action and $s$ is the generator of the equivariant ring of $\C^\times$. 
              \end{thm}
              
	For a toric surface $S$, let $\Delta(S)$ be its Newton polygon and  $\T$  the torus acting on the surface.
	The action of the torus $\T= (\C)^2$ on a toric surface $S$ induces an action of $\T$ on the moduli space of stable pairs. Using Theorem \ref{thm:virtual localisation}, we can reduce the computation to the $\T$-fixed points of the moduli space. The fixed points are isolated and correspond to pairs $(\vec{\lambda},\vec{d})$ such that $\vec{\lambda}\subset\vec{d}$ (see Section \ref{sec:fixed points as lambda, d pairs}), where $\vec{\lambda} = \{\lambda_\alpha\}$ is a collection of partitions indexed by the vertices of the Newton polygon $\Delta(S)$ such that 
	\begin{equation}
		\sum_\alpha |\lambda_\alpha| = m
	\end{equation}
	and $\vec{d}= (d_{\alpha\beta})$ is a tuple of non-negative integers indexed by the edges of $\Delta(S)$ such that 
	\begin{equation}
		[\vec{d}] \coloneqq \sum_{\alpha\beta} d_{\alpha\beta} [D_{\alpha\beta}] = \beta.
	\end{equation}
	Define the shifted series 
	\begin{equation}
		\tilde{Z}_\beta^S(q) \coloneqq q^{\frac{\beta(\beta+K_S)}{2}} Z^S_\beta(q).
	\end{equation}
	Note that by Theorem \ref{thm:virtual localisation},
	\begin{equation} \label{eq:equivariant and non-equivariant series relation}
		\tilde{Z}^{S}_{\beta} (q) = \tilde{Z}^{S}_{\beta} (q\mid s_1,s_2)\big|_{s_1=s_2=0},
	\end{equation}
	where $s_1$ and $s_2$ denote the equivariant parameters. 
 
	By the virtual localisation formula (Theorem \ref{thm:virtual localisation}), it suffices to consider the fixed point locus of the moduli space and sum over those instead of integrating over the entire moduli space. 
	
	\subsection{Charts}\label{sec:charts}

	The vertices of the Newton polygon $\Delta(S)$ correspond to the $\T$-fixed points of $S$ and the edges of $\Delta(S)$ correspond to the $\T$-invariant curves in $S$. 
	The fixed points of the action, as well as the vertices of $\Delta(S)$ are denoted by $p_\alpha$. 
	The $\T$-invariant curve containing the fixed points $p_\alpha$ and $p_\beta$ and the edge connecting the vertices $p_\alpha$ and $p_\beta$ are both denoted by \begin{equation}
		D_{\alpha\beta} = D_{\beta\alpha} \cong\mathbb{P}^1.
	\end{equation}
	Observe that, as $S$ is a surface, to any fixed point $p_\alpha$ we can associate exactly two $\T$-invariant divisors $D_{\alpha\beta_1}$ and $D_{\alpha\beta_2}$ containing $p_\alpha$.
	
	To every fixed point $p_\alpha$, we can associate the canonical $\T$-invariant affine open chart $U_\alpha \cong\mathbb{A}^2$. We denote 
	\begin{equation}
		U_{\alpha\beta} = U_{\beta\alpha} \coloneqq U_\alpha \cap U_\beta.
	\end{equation}
	Let $p_\alpha, p_\beta, p_\gamma, p_\delta$ be vertices in $\Delta(S)$ with edges $D_{\alpha\beta}, D_{\beta\gamma}$, and $D_{\gamma\delta}$. 
	Let the open neighbourhoods of $p_\beta$ and $p_\gamma$ be given by 
	\begin{equation}
		U_\beta \cong \C^2_{x_{\beta\gamma},x_{\beta\alpha}} \qquad \text{ and } \qquad U_\gamma \cong \C^2_{x_{\gamma\beta},x_{\gamma\delta}}.
	\end{equation}
	As $x_{\gamma\beta} = x_{\beta\gamma}^{-1}$, this implies 
	\begin{equation}
		U_{\beta\gamma} \cong \C^\times_{x_{\beta\gamma}} \times \C.
	\end{equation}
	If $D_{\alpha\beta_1}$ and $D_{\alpha\beta_2}$ are the two edges connected to the vertex $p_\alpha$, the action of $\T$ on $U_\alpha$ can be written as:
	\begin{equation}
		(t_{\alpha\beta_1},t_{\alpha\beta_2}) . (x_{\alpha\beta_1},x_{\alpha\beta_2}) = (t_{\alpha\beta_1}x_{\alpha\beta_1}, t_{\alpha\beta_2}x_{\alpha\beta_2}).
	\end{equation} 	
	
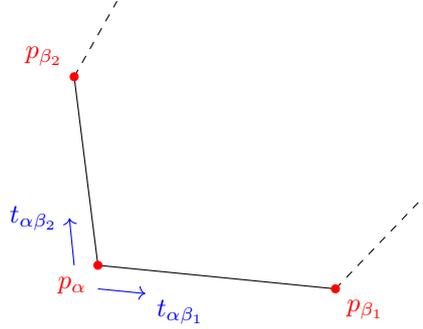
\begin{equation}
		\begin{tikzpicture}[scale=2.5]
			\draw[-] (0,0)--(1.25,-.125);
			\draw[-] (0,0)--(-.125,1);
			\draw[dashed] (1.25,-.125)--(1.75,.4);
			\draw[dashed] (-.125,1)--(.1,1.4);
			\node[below left = 1pt of {(0,0)}] {\textcolor{red}{$p_\alpha$}};
			\node [fill, draw, circle, minimum width=3pt, inner sep=0pt,color=red] at (0,0) {};
			\node[below right = 1pt of {(1.25,-.125)}] {\textcolor{red}{$p_{\beta_1}$}};
			\node [fill, draw, circle, minimum width=3pt, inner sep=0pt,color=red] at (1.25,-.125) {};
			\node[above left= 1pt of {(-.125,1)}] {\textcolor{red}{$p_{\beta_2}$}};
			\node [fill, draw, circle, minimum width=3pt, inner sep=0pt,color=red] at (-.125,1) {};
			\draw[->, color=blue] (0,-.125)--(.25,-.15);
			\node[below right= 1pt of {(.25,-.125)}] {\textcolor{blue}{$t_{\alpha\beta_1}$}};
			\draw[->, color=blue] (-.125,0)--(-.15,.25);
			\node[left = 1pt of {(-.15,.25)}] {\textcolor{blue}{$t_{\alpha\beta_2}$}};
		\end{tikzpicture}
	\end{equation}
	Observe that $t_{\alpha\beta} = (t_{\beta\alpha})^{-1}$.

	\begin{rem}\label{rem: change of variables}
		The normal bundle over a $\T$-invariant curve $D_{\alpha\beta} \cong \mathbb{P}^1$  is given by a line bundle
		\begin{equation}
			\mathcal{N}_{D_{\alpha\beta}/S} =  \mathcal{O}_{\mathbb{P}^1} (m_{\alpha\beta})
		\end{equation} 
		for some $m_{\alpha\beta}\in\mathbb{Z}$.
		Thus, the transition function from $U_\alpha$ to $U_{\beta}$ can be written as 
		\begin{equation}
			(x_{\alpha\beta_1},x_{\alpha\beta_2}) \mapsto (x_{\alpha\beta_1}^{-1}, x_{\alpha\beta_2}x_{\alpha\beta_1}^{-m_{\alpha\beta}}),
		\end{equation} 
		where $\beta_1=\beta$ and $D_{\alpha\beta_2}$ denotes the other edge at $p_\alpha$. Moreover, the curve $D_{\alpha\beta_1}$ is given by $x_{\alpha\beta_2} = 0$. 
	\end{rem}

	\subsection{Monomial Ideals, Young Diagrams, and Partitions}\label{sec:mono,young,part}
	
	\begin{definition}
		A \emph{Young diagram} is a (possibly infinite) collection of boxes arranged in left-justified rows, with the row lengths in non-increasing order. 
	\end{definition}
	
	\begin{ex}
		We consider Young diagrams where the boxes are labelled with Laurent monomials such as 
		\begin{equation}
			\begin{tikzpicture}[scale=1.5]
				\YFrench
				\tyoung(0cm,0cm,1x{x^2}{x^3},y{xy}{x^2y},{y^2}{xy^2},{y^3})
			\end{tikzpicture}
		\end{equation}
		The box labelled $x_1^ix_2^j$ corresponding to $(i,j)$ is the box whose \emph{lower left corner} has coordinates $(i,j)$.
		To each Young diagram we can associate a partition and vice versa. Namely, if a partition includes the tuple $(i,j)$, the Young diagram contains the box labelled $x^iy^j$. 
	\end{ex}
	
	Furthermore, to each Young diagram we can associate the monomial ideal containing all monomials except for those in the diagram. 
	Thus, we have the following correspondences: 
	\begin{equation}
		\begin{tikzcd}
			\bigg\{\T\text{-invariant subschemes}\bigg\} \arrow[r, leftrightarrow, "1:1"] & \bigg\{\text{monomial ideals}\bigg\} \arrow[r, leftrightarrow, "1:1"] & \bigg\{\text{Young diagrams}\bigg\} 
		\end{tikzcd}			
	\end{equation}
	and 
	\begin{equation}
		\begin{tikzcd}
			\bigg\{\begin{array}{c}				0\text{-dimensional }\\				\T\text{-invariant subschemes}			\end{array}\bigg\}\arrow[r, leftrightarrow, "1:1"] & \bigg\{\begin{array}{c}
				\text{finite}\\ \text{Young diagrams}
			\end{array}\bigg\} \arrow[r, leftrightarrow, "1:1"] & \bigg\{\text{partitions}\bigg\}.
		\end{tikzcd}
	\end{equation}
	
	\begin{ex}[Diagram Corresponding to a Divisor $D$]\label{ex:diagram divisor}
		Let $D \subset \mathbb{C}^2$ be a divisor and let $\mathcal{I}_D$ be its monomial ideal given by $(x^{d_{2}}y^{d_{1}})$ where $(d_{2},d_{1})\in\Z^2_{\geq 0}$. The light blue-shaded boxes form the corresponding Young diagram---these extend infinitely upwards and to the right and are referred to as the $(d_{1}, d_{2})$-rays. They are uniquely determined by the coordinate $(d_{2},d_{1})$. 
		\begin{equation}
			\YFrench
			\Yboxdim{1cm}
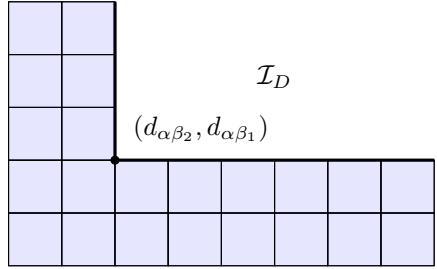
			\begin{tikzpicture}[scale=.7]
				\tyoung(0cm,0cm,~~~~~~~~,~~~~~~~~,~~,~~,~~)
				\draw[-, very thick] (2,2)--(2,5);
				\draw[-, very thick] (2,2)--(8,2);
				\Yfillcolor{blue}
				\Yfillopacity{.1}
				\tyoung(0cm,0cm,~~~~~~~~,~~~~~~~~,~~,~~,~~)
				\node[above right=2pt of {(2,2)}, outer sep=2pt] {$(d_{\alpha\beta_2},d_{\alpha\beta_1})$};
				\node [fill, draw, circle, minimum width=3pt, inner sep=0pt] at (2,2) {};
				\node[above=2pt of {(5,3)}, outer sep=2pt] {$\mathcal{I}_D$};
			\end{tikzpicture}
		\end{equation}
		The partition corresponding to the Young diagram of $D$ is given by
		\begin{equation}
			\{(k_1,k_2) \mid x_1^{k_1}x_2^{k_2} \notin \mathcal{I}_D\} = \{(k_1,k_2) \mid 0\leq k_1 < d_{2} \text{ or } 0\leq k_2 < d_{1}\}.
		\end{equation}
		In particular, $\mathcal{I}_D$ is uniquely determined by the pair $(d_{1},d_{2})$. 
	\end{ex}
	
	\begin{ex}[Diagram Corresponding to a Pair $(Z,D)$]
		Let $D$ be as above in Example \ref{ex:diagram divisor} and $Z$ a $0$-dimensional subscheme contained in $D$. The Young diagram of $Z$ is finite. Note that as $Z\subset D$, we have $\mathcal{I}_D\subset \mathcal{I}_Z$, and hence the Young diagram of $Z$ must be contained in the one of $D$. For example, 
		\begin{equation}
			\YFrench
			\Yboxdim{1cm}
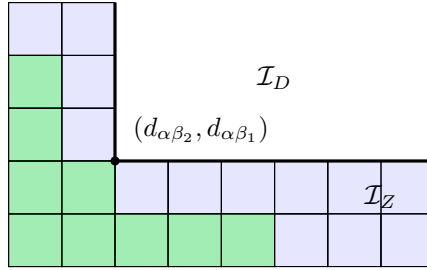
			\begin{tikzpicture}[scale=.7]
				\tyoung(0cm,0cm,~~~~~~~~,~~~~~~,~~,~~,~~)
				\Yfillcolor{blue}
				\Yfillopacity{.1}
				\tyoung(0cm,0cm,~~~~~~~~,~~~~~~~~,~~,~~,~~)
				\Yfillcolor{green}
				\Yfillopacity{.3}
				\tyoung(0cm,0cm,~~~~~,~~,~,~)
				\node[above=2pt of {(7,.75)}, outer sep=2pt] {$\mathcal{I}_Z$};
				\draw[-, very thick] (2,2)--(2,5);
				\draw[-, very thick] (2,2)--(8,2);
				\node [fill, draw, circle, minimum width=3pt, inner sep=0pt] at (2,2) {};
				\node[above right=2pt of {(2,2)}, outer sep=2pt] {$(d_{\alpha\beta_2},d_{\alpha\beta_1})$};
				\node[above=2pt of {(5,3)}, outer sep=2pt] {$\mathcal{I}_D$};
			\end{tikzpicture}
		\end{equation}
		\begin{itemize}
			\item $\mathcal{I}_Z$ contains all monomials outside of the green shaded boxes.
			\item $\mathcal{I}_D$ contains all monomials to the north east of the black lines. Again, the light blue boxes extend infinitely upwards and to the right.
		\end{itemize}
	\end{ex}
	The corresponding partition for $D$ is, as above,
	\begin{equation}
		\{(k_1,k_2) \mid x_1^{k_1}x_2^{k_2} \notin \mathcal{I}_D\} = \{(k_1,k_2) \mid 0\leq k_1 < d_{2} \text{ or } 0\leq k_2 < d_{1}\}
	\end{equation}
	and for $Z$, we have 
	\begin{equation}
		\lambda =  \{(k_1,k_2) \mid x_1^{k_1}x_2^{k_2} \notin \mathcal{I}_Z\}.
	\end{equation}
	
	Note that a pair $(Z,D)$ is uniquely given by the tuple $(d_{1},d_{2})$ and the partition $\lambda$ contained in the $(d_{1},d_{2})$-rays.

	\subsection{Torus Fixed Points of $P_n(S,\beta)$}\label{sec:fixed points as lambda, d pairs}
	
	Let $(Z,D)$ be a pair in the fixed locus of the nested Hilbert scheme under torus action, i.e. $\T$-invariant.  We want to describe this pair in terms of partitions to reduce the computation of $T^\vir|_{(Z\subset D)}$ to a combinatorial problem. Note that Section \ref{sec:mono,young,part} considered this problem over $\mathbb{A}^2$ and we have to \emph{compactify} this discussion, which leads to additional \emph{gluing data} for the $\T$-invariant curves $D_{\alpha\beta}$ corresponding to the edges of the Newton polygon $\Delta(S)$.
	
	Recall that torus-invariant subschemes are given by monomial ideals and vice versa. Hence, a $\T$-invariant divisor $D$ corresponds to a monomial ideal $\mathcal{I}_D$ and a $\T$-invariant $0$-dimensional subscheme $Z$ of $D$ corresponds to a monomial ideal $\mathcal{I}_Z$ containing $\mathcal{I}_D$.
	
	Since $(Z,D)$ is $\T$-fixed on each open set $U_\alpha$, $(\mathcal{I}_Z,\mathcal{I}_D)$ must be given by monomial ideals
	\begin{equation}
		(I_{Z,\alpha},I_{D,\alpha}) =  (\mathcal{I}_Z|_{U_\alpha},\mathcal{I}_D|_{U_\alpha}) \subset\C[x_1,x_2]
	\end{equation} 
	on each $U_\alpha$, too. They may thus (see Section \ref{sec:mono,young,part}), each be viewed as a pair $(d_{\alpha\beta_1},d_{\alpha\beta_2})$ and a partition $\lambda_\alpha$ contained in the $(d_{\alpha\beta_1},d_{\alpha\beta_2})$-rays. 
	
	The integers $d_{\alpha\beta_1}$ and $d_{\alpha\beta_2}$ corresponding to the $\T$-invariant divisor $D_{\alpha\beta_1}$ and $D_{\alpha\beta_2}$, respectively, describes the asymptotic behaviour of the rays in the two coordinate directions.
	This can be viewed as describing the infinite limits of $x_i$-constant cross-sections of the rays of $(d_{\alpha\beta_1},d_{\alpha\beta_2})$. 
	
	To summarise: Every fixed point $(Z,D)$ in $P_n(S,\beta)$ is described fully in terms of
	\begin{enumerate}[label=(\roman*)] 
		\item a tuple $\vec{d} = (d_{\alpha\beta})$, where the $d_{\alpha\beta}$ are non-negative integers and the $\alpha\beta$ run over all edges in the Newton polygon $\Delta(S)$, describing the subscheme $D$ (recall that $D$ is of \emph{pure} dimension $1$) such that 
		\begin{equation}
			\beta= [\vec{d}] \coloneqq \sum_{\alpha\beta} d_{\alpha\beta}[D_{\alpha\beta}],
		\end{equation}
		\item a collection $\vec{\lambda}$ of $2$-dimensional (finite) partitions $\lambda_\alpha$ indexed by $\alpha$ and contained in the $(d_{\alpha\beta_1}, d_{\alpha\beta_2})$-rays, where $p_\alpha$ are the vertices of $\Delta(S)$, corresponding to the subscheme $Z$ such that
		\begin{equation}
			\operatorname{length}(Z) = |\vec{\lambda}| = \sum_\alpha |\lambda_\alpha| = m = n + \frac{\beta(\beta+K_S)}{2}.
		\end{equation}
	\end{enumerate}
	For such $(\vec{\lambda}, \vec{d})$, where the partition $\lambda_\alpha\in\vec{\lambda}$ is contained in the $(d_{\alpha\beta_1},d_{\alpha\beta_2})$-rays, we write $\vec{\lambda}\subset \vec{d}$.

	\section{Virtual Tangent Space}\label{sec:virtual tangent space}
	
	We now want to compute 
	\begin{equation}
		T^\vir \Big|_{Z\subset D} =  \chi(\mathcal{O},\mathcal{I}_Z(D)) - \chi(\mathcal{I}_Z,\mathcal{I}_Z).
              \end{equation}
        
	\subsection{Preliminaries} 
	By \eqref{eq:pot stable pairs} and \eqref{eq:pot nested Hilbert scheme} we see that both perfect obstruction theories have the same virtual tangent space: Let $(F,s)\in P_n(S,\beta)$. Then $\mathcal{O} \xrightarrow{s} F$ can be written as 
	\begin{equation}
		\mathcal{O} \to \mathcal{I}_Z(D) = \mathcal{I}_Z \otimes \mathcal{O}_S(D)
	\end{equation}
	by Proposition \ref{prop:moduli correspondence}. Denote the complex $\{\mathcal{O}_S \to F\}$ by $I^\bullet$ and its derived dual by $(I^\bullet)^\lor$.
	Thus, as already observed in \cite[Proposition 3.1.5]{gsy},
	\begin{align} 
		T^\vir_{P_n(S,\beta)}\Big|_{(F,s)} &= \chi(I^\bullet, F)\\
		&= \chi(I^\bullet, \mathcal{O} - I^\bullet)\\
		&= \chi((\mathcal{O}-I^\bullet)^\lor, (I^\bullet)^\lor)\\
		&= \chi(\mathcal{O} - \mathcal{I}_Z(D), \mathcal{I}_Z(D))\\
		&= T^\vir_{S^{[0,m]}_\beta} \Big|_{Z\subset D},
	\end{align}
	where the second to last equality used that the derived dual is linear and \cite[Lemma A.4]{kt}. 
	Locally, the virtual tangent space at a point $(Z,D) \in P_n(S,\beta)^\T$ corresponding to $(\vec{\lambda}, \vec{d})$ is given by 
	\begin{align}\label{eq:local virtual tangent space}
		T^\vir\Big|_{\vec{\lambda} \subset  \vec{d}} &= \chi(\mathcal{O},\mathcal{I}_{\vec{\lambda}}(\vec{d})) - \chi(\mathcal{I}_{\vec{\lambda}},\mathcal{I}_{\vec{\lambda}})\\
		&= - \chi(\mathcal{O}_{\vec{d}},\mathcal{O}) + \chi(\mathcal{O}_{\vec{d}},\mathcal{O}_{\vec{\lambda}}) + \chi(\mathcal{O}_{\vec{\lambda}},\mathcal{O}) - \chi(\mathcal{O}_{\vec{\lambda}},\mathcal{O}_{\vec{\lambda}}),
	\end{align}
	where we used the exact sequences
	\begin{equation}
		0\to\mathcal{I}_{\vec{\lambda}} \to \mathcal{O} \to \mathcal{O}_{\vec{\lambda}} \to 0
	\end{equation}
	and 
	\begin{equation}
		0 \to \mathcal{O}(-\vec{d}\,) \to \mathcal{O} \to \mathcal{O}_{\vec{d}} \to 0.
	\end{equation}
	We frequently use $\vec{\lambda}$ and $\vec{d}$ is in place of $Z$ and $D$, respectively, to emphasise the partition point of view. 
	
	\subsection{Local-to-Global}\label{sec:local to global}
	We now separate \eqref{eq:local virtual tangent space} into vertex and edge contributions using the local-to-global spectral sequence and \v{C}ech cohomology. For two coherent sheaves $F$ and $G$ on $S$, we have
	\begin{align}\label{eq:local-to-global+cech}
		\chi(F,G) &= H^\bullet(S,\sExt^\bullet (F,))\\
		&=\sum_\alpha H^0(U_\alpha,\sExt^\bullet(F,G)|_{U_\alpha}) - \sum_{\alpha\beta} H^0(U_{\alpha\beta},\sExt^\bullet(F,G)|_{U_{\alpha\beta}})\\
		&= \sum_\alpha H^0(U_\alpha,\sExt^\bullet(F|_{U_\alpha},G|_{U_\alpha})) - \sum_{\alpha\beta} H^0(U_{\alpha\beta}, \sExt^\bullet(F|_{U_{\alpha\beta}},G|_{U_{\alpha\beta}})),
	\end{align}
	where in the first equality we used the local-to-global spectral sequence, in the second equality we used \v{C}ech cohomology and in the last equality we used a fundamental property of the local $\sExt$ sheaf. 
	The vertex contribution  depends on the $\lambda_\alpha$'s and $d_{\alpha\beta}$'s, whereas the edge contribution depends only on the $d_{\alpha\beta}$'s. 
	
	Thus, combining \eqref{eq:local virtual tangent space} and \eqref{eq:local-to-global+cech},
	\begin{multline}
		T^\vir\Big|_{\vec{\lambda} \subset \vec{d}} = \sum_{\alpha} \Big(-H^0(U_\alpha,\sExt^\bullet(\mathcal{O}_{\vec{d}}|_{U_\alpha},\mathcal{O}|_{U_\alpha}))  +  H^0(U_\alpha,\sExt^\bullet(\mathcal{O}_{\vec{d}}|_{U_\alpha},\mathcal{O}_{\vec{\lambda}}|_{U_\alpha})) \\
		+ H^0(U_\alpha,\sExt^\bullet(\mathcal{O}_{\vec{\lambda}}|_{U_\alpha},\mathcal{O}|_{U_\alpha})  - H^0(U_\alpha,\sExt^\bullet(\mathcal{O}_{\vec{\lambda}}|_{U_\alpha},\mathcal{O}_{\vec{\lambda}}|_{U_\alpha}))\Big) \\
		+ \sum_{\alpha\beta} H^0(U_{\alpha\beta},\sExt^\bullet(\mathcal{O}_{\vec{d}}|_{U_{\alpha\beta}},\mathcal{O}|_{U_{\alpha\beta}})).
	\end{multline}
	For every fixed point $p_\alpha$, define the vertex contribution to be
	$$\tilde{\V}_{\lambda_\alpha,d_{\alpha\beta_1}, d_{\alpha\beta_2} } \coloneqq \tilde{\mathsf{F}}_{d_{\alpha\beta_1}, d_{\alpha\beta_2}}+ \mathsf{G}_{\lambda_\alpha}, $$
	where 
	\begin{equation}
		\tilde{\mathsf{F}}_{d_{\alpha\beta_1}, d_{\alpha\beta_2}} \coloneqq -H^0(U_\alpha,\sExt^\bullet(\mathcal{O}_{\vec{d}}|_{U_\alpha},\mathcal{O}|_{U_\alpha}))  
	\end{equation}
	and 
	\begin{multline}
		\mathsf{G}_{\lambda_\alpha}\coloneqq  H^0(U_\alpha,\sExt^\bullet(\mathcal{O}_{\vec{d}}|_{U_\alpha},\mathcal{O}_{\vec{\lambda}}|_{U_\alpha})) + H^0(U_\alpha,\sExt^\bullet(\mathcal{O}_{\vec{\lambda}}|_{U_\alpha},\mathcal{O}|_{U_\alpha})  \\
		-H^0(U_\alpha,\sExt^\bullet(\mathcal{O}_{\vec{\lambda}}|_{U_\alpha},\mathcal{O}_{\vec{\lambda}}|_{U_\alpha})).
	\end{multline}
	For every edge $D_{\alpha\beta}$, define the edge contribution to be
	\begin{equation}
		\tilde{\mathsf{E}}_{d_{\alpha\beta}} \coloneqq H^0(U_{\alpha\beta},\sExt^\bullet(\mathcal{O}_{\vec{d}}|_{U_{\alpha\beta}},\mathcal{O}|_{U_{\alpha\beta}})),
	\end{equation}
	so that 
	\begin{equation}
		T^\vir\Big|_{\vec{\lambda}\subset \vec{d}} = \sum_{\alpha} \tilde{\V}_\alpha + \sum_{\alpha\beta} \tilde{\mathsf{E}}_{\alpha\beta}.
	\end{equation}
	If $\vec{\lambda}$ and $\vec{d}$ are clear from context, we omit them in the notation: 
	\begin{equation}
		\tilde{\V}_\alpha = \tilde{\V}_{\lambda_\alpha}, \tilde{\mathsf{F}}_{d_{\alpha\beta_1}, d_{\alpha\beta_2}} = \tilde{\mathsf{F}}_\alpha, \text{ etc.}
	\end{equation}

	\subsection{Computing $\chi(\mathcal{O}_{\vec{d}}|_{U_\alpha}), \chi(\mathcal{O}_{\vec{d}}|_{U_{\alpha\beta}}), $ and $\chi(\mathcal{O}_{\vec{\lambda}}|_{U_\alpha})$}\label{sec:computations of fundamental ext-terms}

	\subsubsection{Computing $\chi(\mathcal{O}_{\vec{d}}|_{U_{\alpha}})$}\label{sec:Od Ualpha}
	
	Let 
	\begin{equation}
		d_\alpha \coloneqq \{(k_1, k_2) \in\Z_{\geq0}^2 \mid 0 \leq k_1 < d_{\alpha\beta_2} \text{ or } 0\leq k_2 < d_{\alpha\beta_1} \}.
	\end{equation}
	Then
	\begin{equation}
		\chi(\mathcal{O}_{\vec{d}}|_{U_\alpha}) = Q_{d_{\alpha}} = \sum_{(i,j)\in d_\alpha} t_{\alpha\beta_1}^i t_{\alpha\beta_2}^j = Q_{(0,d_{\alpha\beta_1})} + Q_{(d_{\alpha\beta_2},0)} - \sum_{k=1}^{d_{\alpha\beta_2}} \sum_{l=1}^{d_{\alpha\beta_1}} t_{\alpha\beta_1}^{k} t_{\alpha\beta_2}^{l},
	\end{equation}
	where 
	\begin{equation}
		Q_{(0,d_{\alpha\beta_1})} = \sum_{k\in\Z_{\geq0}} t_{\alpha\beta_1}^k \sum_{l=0}^{d_{\alpha\beta_1}-1} t_{\alpha\beta_2}^l.
	\end{equation}
	and
	\begin{equation}
		Q_{(d_{\alpha\beta_2},0)} = \sum_{k=0}^{d_{\alpha\beta_2}-1} t_{\alpha\beta_1}^k \sum_{l\in\Z_{\geq0}} t_{\alpha\beta_2}^l.
	\end{equation}
	Thus,
	\begin{align}
		Q_{d_\alpha} &= \sum_{k\in\Z_{\geq0}} t_{\alpha\beta_1}^k \sum_{l=0}^{d_{\alpha\beta_1}-1} t_{\alpha\beta_2}^l  +  \sum_{k=0}^{d_{\alpha\beta_2}-1} t_{\alpha\beta_1}^k \sum_{l\in\Z_{\geq0}} t_{\alpha\beta_2}^l - \sum_{k=1}^{d_{\alpha\beta_2}} \sum_{l=1}^{d_{\alpha\beta_1}} t_{\alpha\beta_1}^k t_{\alpha\beta_2}^l\\
		&= (1-t_{\alpha\beta_1}^{d_{\alpha\beta_2}}t_{\alpha\beta_2}^{d_{\alpha\beta_1}}) \chi(\mathcal{O}_{U_\alpha}).
	\end{align}

	\subsubsection{Computing $\chi(\mathcal{O}_{\vec{d}}|_{U_{\alpha\beta}})$}

	Recall that to every vertex $\alpha$, we have two associated edges corresponding to the curves $D_{\alpha\beta_1}$ and $D_{\alpha\beta_2}$ and their thickenings are given by $d_{\alpha\beta_1}$ and $d_{\alpha\beta_2}$.  
	Therefore,
	\begin{equation}
		\chi(\mathcal{O}_{\vec{d}}|_{U_{\alpha\beta_1}}) = \sum_{k\in\Z} t_{\alpha\beta_1}^k  \sum_{l=0}^{d_{\alpha\beta_1}-1} t_{\alpha\beta_2}^{l} \quad \text{ and } \quad \chi(\mathcal{O}_{\vec{d}}|_{U_{\alpha\beta_2}}) =  \sum_{l=0}^{d_{\alpha\beta_2}-1} t_{\alpha\beta_1}^{l} \sum_{k\in\Z} t_{\alpha\beta_2}^k .
	\end{equation}
	Let
	\begin{equation}
		\delta(t_{\alpha\beta_i}) \coloneqq \sum_{k\in\Z} t_{\alpha\beta_i}^k.
	\end{equation}
	Then 
	\begin{equation}
		\chi(\mathcal{O}_{\vec{d}}|_{U_{\alpha\beta_i}}) = \delta(t_{\alpha\beta_i}) \left( \sum_{l=0}^{d_{\alpha\beta_i}-1} t_{\alpha\beta_j}^{l} \right) 
	\end{equation}
	for $\{i,j\} = \{1,2\}$.

	\subsubsection{Computing $\chi(\mathcal{O}_{\vec{\lambda}}|_{U_\alpha})$}
	
	By the restriction of $\vec{\lambda}$ to $U_\alpha$ we mean the partition $\lambda_\alpha$, which is a component of  $\vec{\lambda}$. 
	
	The $\T$-character of $\chi(\mathcal{O}_{\lambda_\alpha}) = \chi(\mathcal{O}_{\vec{\lambda}}|_{U_\alpha})$ is given by
	\begin{equation}\label{eq:Vpi_alpha} 
		\chi(\mathcal{O}_{\lambda_\alpha}) = \sum_{(k_1,k_2) \in \lambda_\alpha} t_{\alpha\beta_1}^{k_1} t_{\alpha\beta_2}^{k_2}.
	\end{equation}

	\subsection{Vertex Calculation}
	
	Let $U_\alpha$ be a neighbourhood of a fixed point $p_\alpha$. Using the results of Section  \ref{sec:computations of fundamental ext-terms}, we compute 
	\begin{equation}
		\chi(\mathcal{O}_{\vec{d}}|_{U_\alpha},\mathcal{O}) = \frac{\chi(\mathcal{O}_{\vec{d}}|_{U_\alpha},\mathcal{O})^* \chi(\mathcal{O})}{\chi(\mathcal{O})^*} = \frac{\overline{Q}_{d_\alpha}}{t_{\alpha\beta_1}t_{\alpha\beta_2}}
	\end{equation}
	We have
	\begin{align}
		\tilde{\mathsf{F}}_\alpha &= -\frac{\chi(\mathcal{O}_{\vec{d}}|_{U_\alpha},\mathcal{O})^* \chi(\mathcal{O}_{U_\alpha})}{\chi(\mathcal{O}_{U_\alpha})^*}\\ 
		&= - t_{\alpha\beta_1}^{-1} \sum_{k\geq0}  t_{\alpha\beta_1}^{-k}  \sum_{l=1}^{d_{\alpha\beta_1}} t_{\alpha\beta_2}^{-l}  - t_{\alpha\beta_2}^{-1} \sum_{k\geq0} t_{\alpha\beta_2}^{-k} \sum_{l=1}^{d_{\alpha\beta_2}} t_{\alpha\beta_1}^{-l}  + \sum_{k=1}^{d_{\alpha\beta_2}} \sum_{l=1}^{d_{\alpha\beta_1}} t_{\alpha\beta_1}^k t_{\alpha\beta_2}^l\\
		&= - t_{\alpha\beta_1}^{-1} \frac{1}{1-t_{\alpha\beta_1}^{-1}} \sum_{l=1}^{d_{\alpha\beta_1}} t_{\alpha\beta_2}^{-l}  - t_{\alpha\beta_2}^{-1} \frac{1}{1-t_{\alpha\beta_2}^{-1}} \sum_{l=1}^{d_{\alpha\beta_2}} t_{\alpha\beta_1}^{-l}  + \sum_{k=1}^{d_{\alpha\beta_2}} \sum_{l=1}^{d_{\alpha\beta_1}} t_{\alpha\beta_1}^k t_{\alpha\beta_2}^l\\
		&= - \frac{(1- t_{\alpha\beta_1}^{-d_{\alpha\beta_2}})(1- t_{\alpha\beta_2}^{-d_{\alpha\beta_1}})}{(1- t_{\alpha\beta_1})(1- t_{\alpha\beta_2})}.
	\end{align}
	Furthermore,
	\begin{equation}
		\chi(\mathcal{O}_{\vec{d}}|_{U_\alpha},\mathcal{O}_{\lambda_\alpha}) = (1-t_{\alpha\beta_1}^{-d_{\alpha\beta_2}}t_{\alpha\beta_2}^{-d_{\alpha\beta_1}}) Q_{{\lambda_\alpha}}, 
	\end{equation}
	\begin{equation}
		\chi(\mathcal{O}_{\lambda_\alpha},\mathcal{O}_{U_\alpha}) =\frac{ \overline{Q}_{\lambda_\alpha}}{t_{\alpha\beta_1}t_{\alpha\beta_2}}, 
	\end{equation}
	and
	\begin{equation}
		\chi(\mathcal{O}_{\lambda_\alpha},\mathcal{O}_{\lambda_\alpha}) = (1-t_{\alpha\beta_1}^{-1})(1-t_{\alpha\beta_2}^{-1}) Q_{\lambda_\alpha}\overline{Q}_{\lambda_\alpha} = \frac{(1-t_{\alpha\beta_1})(1-t_{\alpha\beta_2})}{t_{\alpha\beta_1}t_{\alpha\beta_2}} Q_\lambda \overline{Q}_\lambda.  
	\end{equation}
	Thus,
	\begin{equation}\label{eq:Galpha}
		\mathsf{G}_\alpha= (1-t_{\alpha\beta_1}^{-d_{\alpha\beta_2}}t_{\alpha\beta_2}^{-d_{\alpha\beta_1}}) Q_{\lambda_\alpha} + \frac{\overline{Q}_{\lambda_\alpha}}{t_{\alpha\beta_1}t_{\alpha\beta_2}} - \frac{(1-t_{\alpha\beta_1})(1-t_{\alpha\beta_2})}{t_{\alpha\beta_1}t_{\alpha\beta_2}} Q_{\lambda_\alpha} \overline{Q}_{\lambda_\alpha}. 
	\end{equation}
	Whence,
	\begin{multline}\label{eq:V tilde}
		\tilde{\mathsf{F}}_\alpha + \mathsf{G}_\alpha =- t_{\alpha\beta_1}^{-1} \frac{1}{1-t_{\alpha\beta_1}^{-1}} \sum_{l=1}^{d_{\alpha\beta_1}} t_{\alpha\beta_2}^{-l}  - t_{\alpha\beta_2}^{-1} \frac{1}{1-t_{\alpha\beta_2}^{-1}} \sum_{l=1}^{d_{\alpha\beta_2}} t_{\alpha\beta_1}^{-l}  + \sum_{k=1}^{d_{\alpha\beta_2}} \sum_{l=1}^{d_{\alpha\beta_1}} t_{\alpha\beta_1}^k t_{\alpha\beta_2}^l\\ 
		+ (1-t_{\alpha\beta_1}^{-d_{\alpha\beta_2}}t_{\alpha\beta_2}^{-d_{\alpha\beta_1}}) Q_{\lambda_\alpha} + \frac{\overline{Q}_{\lambda_\alpha}}{t_{\alpha\beta_1}t_{\alpha\beta_2}} - \frac{(1-t_{\alpha\beta_1})(1-t_{\alpha\beta_2})}{t_{\alpha\beta_1}t_{\alpha\beta_2}} Q_{\lambda_\alpha} \overline{Q}_{\lambda_\alpha}. 
	\end{multline}

	\subsection{Edge Calculation}
	
	Let $U_{\alpha\beta_i} \coloneqq U_\alpha \cap U_{\beta_i}$, as above. We want to determine the $\T$-character of $\chi(\mathcal{O}|_{U_{\alpha\beta_i}})$. As $U_\alpha\cong\C^2$, the $(k_1,k_2)$ are summed over all of $\Z^2_{\geq0}$ and 
	\begin{equation}
		\chi(\mathcal{O}_{U_\alpha}) = \sum_{(k_1,k_2)\in\Z_{\geq0}} t_{\alpha\beta_1}^{k_1} t_{\alpha\beta_2}^{k_2}.
              \end{equation}
        Recall that if $U_\alpha \cong \C_{x_{\alpha\beta_i}} \times \C$ and 
	$U_{\beta_i} \cong \C_{x_{\alpha\beta_i}^{-1}} \times \C$, then $U_{\alpha\beta_i} \cong \C_{x_{\alpha\beta_i}}^\times \times \C$ and so  $k_i \neq 0$. 
	Therefore, 
	\begin{equation}
		\chi(\mathcal{O}|_{U_{\alpha\beta_i}}) = \delta(t_{\alpha\beta_i}) \sum_{k\geq0}t_{\alpha\beta_j}^k = \frac{\delta(t_{\alpha\beta_i})}{1-t_{\alpha\beta_j}}.
	\end{equation}
	Thus, the edge contribution is 
	\begin{align}
		\mathsf{\tilde{E}}_{\alpha\beta_i}  &= - \chi(\mathcal{O}_{\vec{d}}|_{U_{\alpha\beta_i}},\mathcal{O}|_{U_{\alpha\beta_i}})  \\ 
		&= - \frac{\chi(\mathcal{O}|_{U_{\alpha\beta_i}})}{\chi(\mathcal{O}|_{U_{\alpha\beta_i}})^*} \cdot \delta(t_{\alpha\beta_i}) \left( \sum_{l=0}^{d_{\alpha\beta_i}-1} t_{\alpha\beta_j}^{-l}\right)		\\
		&= - \delta(t_{\alpha\beta_i}) \sum_{l=1}^{d_{\alpha\beta_i}} t_{\alpha\beta_j}^{-l},\label{eq:E tilde}
	\end{align}
	where $\{i,j\}  = \{1,2\}$ and in the last equation we used that
	\begin{equation}
		\frac{\chi(\mathcal{O}|_{U_{\alpha\beta_i}})}{\chi(\mathcal{O}|_{U_{\alpha\beta_i}})^*} = \frac{1}{t_{\alpha\beta_j}}
	\end{equation}
	and
	\begin{equation}
		\delta(t_{\alpha\beta_i}) = t_{\alpha\beta_i}\delta(t_{\alpha\beta_i}) = \delta(t_{\alpha\beta_i}^{-1}).
	\end{equation}
	To ease notation, define 
	\begin{equation}
		\mathsf{F}_{\alpha\beta_i} (t)\coloneqq -\sum_{l=1}^{d_{\alpha\beta_i}} t^{-l}. 
	\end{equation}
	That is, let
	\begin{equation}
		\mathsf{\tilde{E}}_{\alpha\beta_i} = \delta(t_{\alpha\beta_i}) \mathsf{F}_{\alpha\beta_i} (t_{\alpha\beta_j})
	\end{equation}
	
	\begin{rem}
		We have so far computed the $\T$-character of the virtual tangent space as the sum
		\begin{equation}
			T^\vir = \sum_\alpha (\tilde{\mathsf{F}}_\alpha+ \mathsf{G}_\alpha) +\sum_{\alpha\beta} \tilde{\mathsf{E}}_{\alpha\beta}.
		\end{equation}
		
		The problem is that the $\tilde{\mathsf{V}}_\alpha$ and $\tilde{\mathsf{E}}_{\alpha\beta}$ are not Laurent polynomials, but rather formal Laurent series. Note that the problem in $\tilde{\V}_\alpha$ lies only in the term $\tilde{\mathsf{F}}_\alpha$ depending solely on $d_\alpha$ and that $\tilde{\mathsf{E}}_{\alpha\beta}$ depends only on $d_{\alpha\beta}$.  The solution is to use a \emph{change of variables} coming from the normal bundle to regroup the individual terms of $\tilde{\mathsf{F}}_\alpha$ and $\tilde{\mathsf{E}}_{\alpha\beta}$, defining $\mathsf{F}_\alpha$ and $\mathsf{E}_{\alpha\beta}$ such that 
		\begin{equation}
			T^\vir = \sum_\alpha (\mathsf{F}_\alpha + \mathsf{G}_\alpha ) +  \sum_{\alpha\beta} \mathsf{E}_{\alpha\beta},
		\end{equation}
		which are then shown to be Laurent polynomials in Lemma \ref{lemma:Laurent}. 
		
	\end{rem}
	
	\subsection{Regrouping}
	
	\subsubsection{Change of Variables} 
	We regroup  the terms of the vertex and edge contributions to get two finite sums.

	The edge contribution can we written as the sum of two terms, one summing over all non-negative terms, and one over all negative terms:
	\begin{equation}
		\tilde{\mathsf{E}}_{\alpha\beta} = \mathsf{P}_{\alpha\beta} + \mathsf{N}_{\alpha\beta}  = \frac{\mathsf{F}_{\alpha\beta}}{1-t_{\alpha\beta}} + t_{\alpha\beta}^{-1} \frac{\mathsf{F}_{\alpha\beta}}{1-t_{\alpha\beta}^{-1}}.
	\end{equation}
	Note that here the ordering of $\alpha$ and $\beta$ in the indices is relevant, a change of variables will lead to the same result using $\beta\alpha$:
	\begin{equation}
		\tilde{\mathsf{E}}_{\alpha\beta} = \mathsf{P}_{\beta\alpha}+ \mathsf{N}_{\beta\alpha}
	\end{equation}
	The change of variables 
	\begin{equation}
		(x_1,x_2) \mapsto (x_1^{-1},x_2x_1^{-m_{\alpha\beta}})
	\end{equation}
	from Remark \ref{rem: change of variables} relates $\mathsf{P}_{\beta\alpha}$ to $\mathsf{N}_{\alpha\beta}$ and $\mathsf{P}_{\alpha\beta}$ to $\mathsf{N}_{\beta\alpha}$.

	\subsubsection{Regrouping}\label{sec:regrouping}

	\begin{lemma}\label{lemma:d vertex contribution is Laurent}
		The term $\tilde{\mathsf{F}}_\alpha$ can be written as 
		\begin{equation}
			\tilde{\mathsf{F}}_\alpha =  t_{\alpha\beta_1}^{-1}\frac{\mathsf{F}_{\alpha\beta_1}}{1-t_{\alpha\beta_1}^{-1}} + \ldots  
		\end{equation}
		where the dots stand for terms with no pole at $t_{\alpha\beta_1} = 1$. 
	\end{lemma}
	
	\begin{proof}
		Indeed, 
		\begin{align}
			\tilde{\mathsf{F}}_\alpha &= - t_{\alpha\beta_1}^{-1} \frac{1}{1-t_{\alpha\beta_1}^{-1}} \sum_{l=1}^{d_{\alpha\beta_1}} t_{\alpha\beta_2}^{-l}  - t_{\alpha\beta_2}^{-1} \frac{1}{1-t_{\alpha\beta_2}^{-1}} \sum_{l=1}^{d_{\alpha\beta_2}} t_{\alpha\beta_1}^{-l}  +\sum_{k=1}^{d_{\alpha\beta_2}} \sum_{l=1}^{d_{\alpha\beta_1}} t_{\alpha\beta_1}^k t_{\alpha\beta_2}^l\\
			&= t_{\alpha\beta_1}^{-1} \frac{\mathsf{F}_{\alpha\beta_1}(t_{\alpha\beta_2})}{1-t_{\alpha\beta_1}^{-1}} + t_{\alpha\beta_2}^{-1}\frac{\mathsf{F}_{\alpha\beta_2}(t_{\alpha\beta_1})}{1-t_{\alpha\beta_2}^{-1}} + \sum_{j=1}^{d_{\alpha\beta_2}} \sum_{k=1}^{d_{\alpha\beta_1}} t_{\alpha\beta_1}^{j} t_{\alpha\beta_2}^{k}.
		\end{align}
		The latter two terms are clearly regular at $t_{\alpha\beta_1} =1$. 
	\end{proof}
	Define 
	\begin{equation}\label{eq:vertex contribution}
		\V_\alpha \coloneqq \mathsf{F}_\alpha +  \mathsf{G}_\alpha - \sum_{i=1}^2 t_{\alpha\beta_i}^{-1}\frac{\mathsf{F}_{\alpha\beta_i}(t_{\alpha\beta_j})}{1-t_{\alpha\beta_i}^{-1}} =  \mathsf{G}_\alpha +\sum_{k=1}^{d_{\alpha\beta_2}} \sum_{l=1}^{d_{\alpha\beta_1}} t_{\alpha\beta_1}^k t_{\alpha\beta_2}^l
	\end{equation}
	and 
	\begin{equation}\label{eq:edge contribution}
		\mathsf{E}_{\alpha\beta_i} \coloneqq t_{\alpha\beta_i}^{-1} \frac{\mathsf{F}_{\alpha\beta_i}(t_{\alpha\beta_j})}{1-t_{\alpha\beta_i}^{-1}} - \frac{\mathsf{F}_{\alpha\beta_i}(t_{\alpha\beta_j}t_{\alpha\beta_i}^{-m_{\alpha\beta_i}})}{1-t_{\alpha\beta_i}^{-1}}. 
	\end{equation}
	Note that the $\mathsf{E}_{\alpha\beta}$ are  symmetric.
	
	\begin{lemma}\label{lemma:Laurent}
		$\V_\alpha$ and $\mathsf{E}_{\alpha\beta_i}$ are Laurent polynomials.
	\end{lemma}
	
	\begin{proof}
		Clearly, $\mathsf{G}_\alpha$ is a Laurent polynomial. By Lemma \ref{lemma:d vertex contribution is Laurent},
		\begin{equation}
			\V_\alpha = \mathsf{G}_\alpha - t_{\alpha\beta_2}^{-1} \frac{\mathsf{F}_{\alpha\beta_2}}{1-t_{\alpha\beta_2}^{-1}} + \ldots,
		\end{equation}
		where the dots stand for terms regular at $t_{\alpha\beta_1} =1$. 
		Hence, $\V_\alpha$ is a Laurent polynomial. Furthermore, note that the numerator of $\mathsf{E}_{\alpha\beta_1}$ vanishes at $t_{\alpha\beta_1} = 1$. 
	\end{proof}
	
	\begin{thm}\label{thm:virtual tangent space}
		The $\T$-character of $T^\vir$ is given by 
		\begin{equation}
			T^\vir (t_{1},t_{2}) = \sum_\alpha \V_\alpha + \sum_{\alpha\beta} \mathsf{E}_{\alpha\beta}.
		\end{equation}
	\end{thm}
	
	\begin{proof}
		This follows from the definition of $\V_\alpha$ and $\mathsf{E}_{\alpha\beta}$ and Lemma \ref{lemma:Laurent}.
	\end{proof}

	\begin{rem} 
		Note that $\V_\alpha$ is associated to the Hilbert scheme of the singular curve $C$ associated to $\vec{d}$.  
		The fixed points of $C_{\vec{d}}^{[m]}$ are given by 
		\begin{equation}
			\left( C^{[m]}_{\vec{d}} \right)^\T\simeq  \{\vec{\lambda} = \{\lambda_\alpha\} \subset \vec{d}\}.
		\end{equation}
		Thus, if we write the virtual tangent space of $C_{\vec{d}}^{[n]}$ as 
		\begin{equation}
			T^\vir (S^{[m]}) \Big|_{\vec{\lambda}} - L_\beta^{[m]}\Big|_{\vec{\lambda}},
		\end{equation}
		the first term is given by 
		\begin{equation}
			\sum_\alpha \left(Q_\alpha + t_{\alpha\beta_1}^{-1} t_{\alpha\beta_2}^{-1} \overline{Q}_\alpha - (1-t_{\alpha\beta_1}^{-1})(1-t_{\alpha\beta_2}^{-1}) Q_\alpha\overline{Q}_\alpha\right)
		\end{equation}
		and the second term by
		\begin{align}
			\chi(\mathcal{O}(\vec{d}) \otimes \mathcal{O}_{\vec{\lambda}}) &= \chi(\mathcal{O}(-\vec{d}), \mathcal{O}_{\vec{\lambda}}) \\
			&= \sum_\alpha H^0(U_\alpha, \sExt^\bullet(\mathcal{O}(-\vec{d}), \mathcal{O}_{\vec{\lambda}}))\\&= \sum_\alpha Q_\alpha t_{\alpha\beta_1}^{-d_{\alpha\beta_2}} t_{\alpha\beta_2}^{-d_{\alpha\beta_1}}.
		\end{align}
	\end{rem}

	\subsection{Vanishing Fixed Part of the Virtual Tangent Space}
	
	\begin{thm}\label{thm:Tvir no fixed part}
		The virtual tangent space $T^\vir$ does not have a fixed part.
	\end{thm}
	
	\begin{proof}
		By Theorem \ref{thm:virtual tangent space},
		\begin{equation}
			T^\vir (t_{1},t_{2}) = \sum_\alpha \V_\alpha + \sum_{\alpha\beta} \mathsf{E}_{\alpha\beta}.
		\end{equation}
		The edge contribution does not have a fixed part, as the sums $\mathsf{F}_{\alpha\beta}$ do not have a $t_{\alpha\beta_1}^0t_{\alpha\beta_2}^0$ term. 
		Note that ${\V_\alpha}$ is given by
		\begin{equation}
			\sum_{k=1}^{d_{\alpha\beta_2}} \sum_{l=1}^{d_{\alpha\beta_1}} t_{\alpha\beta_1}^k t_{\alpha\beta_2}^l + \underbrace{(1-t_{\alpha\beta_1}^{-d_{\alpha\beta_2}}t_{\alpha\beta_2}^{-d_{\alpha\beta_1}}) Q_{\lambda_\alpha} + \frac{\overline{Q}_{\lambda_\alpha}}{t_{\alpha\beta_1}t_{\alpha\beta_2}} - \frac{(1-t_{\alpha\beta_1})(1-t_{\alpha\beta_2})}{t_{\alpha\beta_1}t_{\alpha\beta_2}} Q_{\lambda_\alpha} \overline{Q}_{\lambda_\alpha}}_{\eqcolon \mathsf{G}_{\alpha}}.
		\end{equation}
		The first term clearly does not have a fixed summand. 
		It remains to show that the fixed part of $\mathsf{G}_\alpha$ vanishes. To this end, we determine the fixed parts of the individual terms of $\mathsf{G}_\alpha$.
		
		First, let us compute $(Q_{\lambda_\alpha})^{\fix}$:
		\begin{equation}
			Q_{\lambda_\alpha} = \sum_{(k_1,k_2) \in \lambda_\alpha} t_{\alpha\beta_1}^{k_1}t_{\alpha\beta_2}^{k_2} = t_{\alpha\beta_1}^0t_{\alpha\beta_2}^0 + \ldots,
		\end{equation}
		where the dots stand for terms dependent on $t_{\alpha\beta_1},t_{\alpha\beta_2}$. Hence,
		\begin{equation}
			(Q_{\lambda_\alpha})^{\fix} = 1.
		\end{equation}
		Furthermore, we determine $(t_{\alpha\beta_1}^{-d_{\alpha\beta_2}}t_{\alpha\beta_2}^{-d_{\alpha\beta_1}}Q_{\lambda_\alpha})^{\fix}$:
		\begin{equation}
			(t_{\alpha\beta_1}^{-d_{\alpha\beta_2}}t_{\alpha\beta_2}^{-d_{\alpha\beta_1}}Q_{\lambda_\alpha})^{\fix} = 
			\left\{ \begin{array}{ll}
				1 & \text{ if } (d_{\alpha\beta_2},d_{\alpha\beta_1}) \in \lambda_\alpha,\\
				0 & \text{ if } (d_{\alpha\beta_2},d_{\alpha\beta_1}) \notin \lambda_\alpha.
			\end{array} \right.
		\end{equation}
		As $(d_{\alpha\beta_2},d_{\alpha\beta_1})$ cannot be in $\lambda_\alpha$, this term does not have a fixed part. 
		Therefore, 
		\begin{equation}\label{eq:fixed part vertex second term}
			((1-t_{\alpha\beta_1}^{-d_{\alpha\beta_2}}t_{\alpha\beta_2}^{-d_{\alpha\beta_1}}) Q_{\lambda_\alpha})^{\fix} = 1.
		\end{equation}
		Second, consider $\frac{\overline{Q}_{\lambda_\alpha}}{t_{\alpha\beta_1}t_{\alpha\beta_2}}$:
		\begin{equation}\label{eq:fixed part vertex third term}
			\frac{\overline{Q}_{\lambda_\alpha}}{t_{\alpha\beta_1}t_{\alpha\beta_2}} = t_{\alpha\beta_1}^{-1}t_{\alpha\beta_2}^{-1} \sum_{(k_1,k_2)\in\lambda_\alpha} t_{\alpha\beta_1}^{-k_1} t_{\alpha\beta_2}^{k_2} = \sum_{(k_1,k_2)\in\lambda_\alpha} t_{\alpha\beta_1}^{-(k_1+1)}t_{\alpha\beta_2}^{-(k_2+1)}.
		\end{equation}
		Hence, this term does not have a fixed part. 
		
		Finally, consider the term 
		\begin{align}
			\frac{(1-t_{\alpha\beta_1})(1-t_{\alpha\beta_2})}{t_{\alpha\beta_1}t_{\alpha\beta_2}} Q_{\lambda_\alpha} \overline{Q}_{\lambda_\alpha} &= (1-t_{\alpha\beta_1}^{-1})(1-t_{\alpha\beta_2}^{-1})  Q_{\lambda_\alpha} \overline{Q}_{\lambda_\alpha}\\
			&= Q_{\lambda_\alpha} \overline{Q}_{\lambda_\alpha} - t_{\alpha\beta_1}^{-1}Q_{\lambda_\alpha} \overline{Q}_{\lambda_\alpha} - t_{\alpha\beta_2}^{-1}Q_{\lambda_\alpha} \overline{Q}_{\lambda_\alpha} + t_{\alpha\beta_1}^{-1}t_{\alpha\beta_2}^{-1}Q_{\lambda_\alpha} \overline{Q}_{\lambda_\alpha}. \qquad \label{eq:Galpha fixed part}
		\end{align}
		Define 
		\begin{equation}
			A_1 \coloneqq \{k_1 \mid \exists k_2\in\Z_{\geq0} :(k_1,k_2) \in \lambda_\alpha\}, \quad A_2 \coloneqq \{k_2 \mid \exists k_1\in\Z_{\geq0} :(k_1,k_2) \in \lambda_\alpha\}.
		\end{equation}
		We compute the fixed part of the individual terms in \eqref{eq:Galpha fixed part}:
		\begin{equation}
			Q_{\lambda_\alpha} \overline{Q}_{\lambda_\alpha} = \sum_{\substack{(k_1,k_2) \in \lambda_\alpha\\(l_1,l_2) \in \lambda_\alpha}} t_{\alpha\beta_1}^{k_1-l_1} t_{\alpha\beta_2}^{k_2-l_2}
		\end{equation}
		The fixed part is where $k_1 = l_1$ and $k_2=l_2$. Hence,
		\begin{equation}\label{eq:fixed part double}
			(Q_{\lambda_\alpha} \overline{Q}_{\lambda_\alpha})^{\fix} = |\lambda_\alpha|.
		\end{equation}
		For the term
		\begin{equation}
			t_{\alpha\beta_1}^{-1}Q_{\lambda_\alpha} \overline{Q}_{\lambda_\alpha}  = t_{\alpha\beta_1}^{-1} \sum_{\substack{(k_1,k_2) \in \lambda_\alpha\\(l_1,l_2) \in \lambda_\alpha}} t_{\alpha\beta_1}^{k_1-l_1} t_{\alpha\beta_2}^{k_2-l_2} = \sum_{\substack{(k_1,k_2) \in \lambda_\alpha\\(l_1,l_2) \in \lambda_\alpha}} t_{\alpha\beta_1}^{k_1-l_1-1} t_{\alpha\beta_2}^{k_2-l_2},
		\end{equation}
		we have that the fixed part is where $l_1 = k_1 - 1$ and $l_2 = k_2$. In particular, as $l_1 \geq 0$, we have $k_1\geq 1$. Therefore, 
		\begin{equation}\label{eq:fixed part double one}
			(t_{\alpha\beta_1}^{-1}Q_{\lambda_\alpha} \overline{Q}_{\lambda_\alpha})^{\fix} = |\lambda_\alpha| - |A_2|.
		\end{equation}
		Analogously, we have
		\begin{equation}\label{eq:fixed part double two}
			(t_{\alpha\beta_2}^{-1}Q_{\lambda_\alpha} \overline{Q}_{\lambda_\alpha})^{\fix} = |\lambda_\alpha| - |A_1|.
		\end{equation}
		Similarly, 
		\begin{equation}\label{eq:fixed part double both}
			(t_{\alpha\beta_1}^{-1}t_{\alpha\beta_2}^{-1}Q_{\lambda_\alpha} \overline{Q}_{\lambda_\alpha})^{\fix} = |\lambda_\alpha| - (|A_1|+|A_2|-1) = |\lambda_\alpha| - |A_1| - |A_2| + 1,
		\end{equation}
		where we deduct $1$ from $|A_1| + |A_2|$, because we count one appearance of $t_{\alpha\beta_1}^0t_{\alpha\beta_2}^0$ twice, once via $A_1$ and once via $A_2$. 
		
		Then, \eqref{eq:fixed part double}, \eqref{eq:fixed part double one}, \eqref{eq:fixed part double two}, and \eqref{eq:fixed part double both} give 
		\begin{equation}\label{eq:fixed part vertex fourth term}
			\left(\frac{(1-t_{\alpha\beta_1})(1-t_{\alpha\beta_2})}{t_{\alpha\beta_1}t_{\alpha\beta_2}} Q_{\lambda_\alpha} \overline{Q}_{\lambda_\alpha}\right)^{\fix} = 1.
		\end{equation}
		Hence, as the first term does not have a fixed part and by \eqref{eq:fixed part vertex second term}, \eqref{eq:fixed part vertex third term}, \eqref{eq:fixed part vertex fourth term}, the fixed part of the virtual tangent space is given by 
		\begin{equation}
			(\V_\alpha)^{\fix} = 0+ 1 + 0 - 1 = 0.
		\end{equation}
		As the edge contribution does not have a fixed part, this concludes the proof. 
	\end{proof}

 \subsection{Proof of Theorem \ref{prop:main conjecture is equiv to loc}}

We now accumulated all the tools necessary to prove Theorem \ref{prop:main conjecture is equiv to loc}. 
 \begin{thm}[= Theorem \ref{prop:main conjecture is equiv to loc}]
	For a non-singular toric surface $S$ and  a homology class $\beta\in H_2(S;\Z)$, 
	\begin{equation}
		\tilde{Z}^{S}_{\beta} (q\mid s_1,s_2)  = \sum_{m\geq 0} \sum_{[\vec{d}] = \beta}   \sum_{\substack{\vec{\lambda}\subset \vec{d} \\|\vec{\lambda}| = m}} q^{m}  \prod_{\alpha\beta} \frac{c_\bullet}{\eeuler}\left(\mathsf{E}_{\alpha\beta}\right) \prod_{\alpha} \frac{c_\bullet}{\eeuler} \left(\mathsf{F}_\alpha + \mathsf{G}_\alpha\right),
	\end{equation}
	where 
	\begin{equation}
		T^\vir = \sum_{\alpha\beta} \mathsf{E}_{\alpha\beta} + \sum_\alpha (\mathsf{F}_\alpha + \mathsf{G}_\alpha)
	\end{equation}
	and $\mathsf{G}_\alpha, \mathsf{F}_\alpha$, and $\mathsf{E}_{\alpha\beta}$ are defined as in \eqref{eq:Galpha}, \eqref{eq:vertex contribution}, and \eqref{eq:edge contribution}, respectively.
\end{thm}

\begin{proof}
    By Theorems \ref{thm:virtual tangent space} and \ref{thm:Tvir no fixed part}, 
    \begin{equation}
		T^\vir = \sum_{\alpha\beta} \mathsf{E}_{\alpha\beta} + \sum_\alpha (\mathsf{F}_\alpha + \mathsf{G}_\alpha).
	\end{equation}
    By Lemma \ref{lemma:Laurent}, this is well-defined. The statement then follows from the definition of $\tilde{Z}_\beta^S(q\mid s_1,s_2)$ given by \eqref{eq:Z tilde},
    \begin{align}
		\tilde{Z}_\beta^S(q\mid s_1,s_2) &= \sum_{m\geq 0} q^m \int_{\left[S_\beta^{[0,m]}\right]^\vir} c_\bullet(T^\vir) \in\mathbb{Q} (s_1,s_2)[[q]]\\
		&= \sum_{[\vec{d}] =\beta} \sum_{\vec{\lambda}\subset \vec{d}} q^{|\vec{\lambda}|} \frac{c_\bullet}{\eeuler} \left(T^\vir\Big|_{\vec{\lambda}\subset\vec{d}}\right)\\
		&= \sum_{[\vec{d}] =\beta} \prod_{\alpha\beta} \frac{c_\bullet}{\eeuler} (\mathsf{E}_{\alpha\beta}) \prod_\alpha   \frac{c_\bullet}{\eeuler} (\mathsf{F}_\alpha) \sum_{\vec{\lambda}\subset\vec{d}} q^{|\vec{\lambda}|} \prod_\alpha \frac{c_\bullet}{\eeuler} (\mathsf{G}_\alpha)\\
        &= \sum_{m\geq 0} \sum_{[\vec{d}] = \beta}   \sum_{\substack{\vec{\lambda}\subset \vec{d} \\|\vec{\lambda}| = m}} q^{m}  \prod_{\alpha\beta} \frac{c_\bullet}{\eeuler}\left(\mathsf{E}_{\alpha\beta}\right) \prod_{\alpha} \frac{c_\bullet}{\eeuler} \left(\mathsf{F}_\alpha + \mathsf{G}_\alpha\right).
	\end{align}
\end{proof}

	\section{Computing $\frac{c_\bullet}{\eeuler} (T^\vir)$}\label{sec:chern euler}
	
	\paragraph{Chern and Euler Class.}
	
	To ease notation we shall denote 
	\begin{equation}
		f(x) \coloneqq \frac{1+x}{x}.
	\end{equation}
	In particular, 
	\begin{equation}
		\frac{c_\bullet}{\eeuler}(t) = f(s)
	\end{equation}
	if $s$ denotes the class of $t$. 
	
	We divide the computation of the virtual tangent space at a point corresponding to $\vec{\lambda}\subset\vec{d}$ into multiple sub-parts:
	
	\begin{equation}
		\frac{c_\bullet}{\eeuler}\left(T^\vir\Big|_{\vec{\lambda}\subset\vec{d}}\right) = \prod_\alpha \frac{c_\bullet}{\eeuler}(\V_\alpha) \prod_{\alpha,\beta} \frac{c_\bullet}{\eeuler} (\mathsf{E}_{\alpha\beta}),
	\end{equation}
	where 
	\begin{equation}
		\frac{c_\bullet}{\eeuler}(\V_\alpha) = \prod_\alpha \frac{c_\bullet}{\eeuler} (\mathsf{G}_\alpha) \cdot \frac{c_\bullet}{\eeuler} \left(-\sum_{k= 1}^{d_{\alpha\beta_2}} \sum_{l=1}^{d_{\alpha\beta_1}} t_{\alpha\beta_1}^{-k} t_{\alpha\beta_2}^{-l}\right).
	\end{equation}
	Observe that 
	\begin{equation}\label{eq:Falpha}
		\frac{c_\bullet}{\eeuler} \left(\sum_{k= 1}^{d_{\alpha\beta_2}} \sum_{l=1}^{d_{\alpha\beta_1}} t_{\alpha\beta_1}^{-k} t_{\alpha\beta_2}^{-l}\right) = \prod_{k=1}^{d_{\alpha\beta_2}} \prod_{l=1}^{d_{\alpha\beta_1}} f(-k s_{\alpha\beta_1} - l s_{\alpha\beta_2}).
	\end{equation}
	By \eqref{eq:Galpha},
	\begin{equation}
		\frac{c_\bullet}{\eeuler}(\mathsf{G}_\alpha) = \underbrace{\frac{c_\bullet}{\eeuler} \left((1-t_{\alpha\beta_1}^{-d_{\alpha\beta_2}}t_{\alpha\beta_2}^{-d_{\alpha\beta_1}})Q_{\lambda_\alpha}\right)}_{\eqqcolon C_1}\cdot \underbrace{\frac{c_\bullet}{\eeuler}\left(\frac{\overline{Q}_{\lambda_\alpha}}{t_{\alpha\beta_1}t_{\alpha\beta_2}}\right)}_{\eqqcolon C_2}\cdot {\underbrace{\frac{c_\bullet}{\eeuler}\left(\frac{(1-t_{\alpha\beta_1})(1-t_{\alpha\beta_2})}{t_{\alpha\beta_1}t_{\alpha\beta_2}} Q_{\lambda_\alpha} \overline{Q}_{\lambda_\alpha}\right)}_{\eqqcolon C_3}}^{-1}.
	\end{equation}
	
	\paragraph{Computing $C_1$.}
	
	Let $s_i$ denote the class of $t_{\alpha\beta_i}$ in cohomology. 
	As 
	\begin{equation}
		C_1 = \frac{c_\bullet}{\eeuler}\left(Q_{\lambda_\alpha}\right) \cdot \frac{c_\bullet}{\eeuler}\left(t_{\alpha\beta_1}^{-d_{\alpha\beta_2}}t_{\alpha\beta_1}^{-d_{\alpha\beta_2}}Q_{\lambda_\alpha}\right)^{-1}
	\end{equation}
	and by \eqref{eq:Vpi_alpha}
	\begin{equation}
		\frac{c_\bullet}{\eeuler}\left(Q_{\lambda_\alpha}\right)  = \frac{c_\bullet}{\eeuler}\left(\sum_{(k,l)\in\lambda_\alpha} t_{\alpha\beta_1}^kt_{\alpha\beta_2}^l \right)
		= \prod_{(k,l)\in\lambda_\alpha}   f(k s_{\alpha\beta_1} + l s_{\alpha\beta_2})
	\end{equation}
	and similarly
	\begin{equation}
		\frac{c_\bullet}{\eeuler}(t_{\alpha\beta_1}^{-d_{\alpha\beta_1}}t_{\alpha\beta_2}^{-d_{\alpha\beta_2}}Q_{\lambda_\alpha}) = \prod_{(k,l)\in\lambda_\alpha}   f((k-d_{\alpha\beta_1})s_{\alpha\beta_1} + (l-d_{\alpha\beta_2}) s_{\alpha\beta_2}).
	\end{equation}
	Thus,
	\begin{equation}\label{eq:C1}
		C_1 = \prod_{(k,l)\in\lambda_\alpha}   \frac{f(k s_{\alpha\beta_1} + l s_{\alpha\beta_2})}{f((k-d_{\alpha\beta_1})s_{\alpha\beta_1} + (l-d_{\alpha\beta_2}) s_{\alpha\beta_2})}.
	\end{equation}
	
	\paragraph{Computing $C_2$.}
	
	Computing analogously as above, we have,
	\begin{align}\label{eq:C2}
		C_2  &= \frac{c_\bullet}{\eeuler}(t_{\alpha\beta_1}^{-1}t_{\alpha\beta_2}^{-1}\overline{Q}_{\lambda_\alpha}) \\
		&= \frac{c_\bullet}{\eeuler} \left(\sum_{(k,l)\in\lambda_\alpha} t_{\alpha\beta_1}^{-(k+1)}t_{\alpha\beta_2}^{-(l+1)}\right)\\
		&= \prod_{(k,l)\in\lambda_\alpha}    f(-(k+1)s_{\alpha\beta_1}- (l+1)s_{\alpha\beta_2}).
	\end{align}
	
	\paragraph{Computing $C_3$.}
	
	First, note that $Q_{{\lambda_\alpha}}\overline{Q}_{\lambda_\alpha}$:
	\begin{equation}
		Q_{\lambda_\alpha}\overline{Q}_{\lambda_\alpha} = \sum_{(k,l)\in\lambda_\alpha} \sum_{(k',l')\in\lambda_\alpha} t_{\alpha\beta_1}^{k-k'} t_{\alpha\beta_2}^{l-l'}.
	\end{equation}
	Thus,
	\begin{equation}
		\frac{Q_{{\lambda_\alpha}}\overline{Q}_{\lambda_\alpha}}{t_{\alpha\beta_1}t_{\alpha\beta_2}} =  \sum_{(k,l)\in\lambda_\alpha} \sum_{(k',l')\in\lambda_\alpha} t_{\alpha\beta_1}^{k-k'-1} t_{\alpha\beta_2}^{l-l'-1}.
	\end{equation}
	Therefore, the individual terms of which $C_3$ consists are
	\begin{align}
		\frac{c_\bullet}{\eeuler} \left( \frac{Q_{{\lambda_\alpha}}\overline{Q}_{\lambda_\alpha}}{t_{\alpha\beta_1}t_{\alpha\beta_2}}\right) &= \prod_{k,k', l,l'} f((k-k'-1)s_{\alpha\beta_1} + (l-l'-1) s_{\alpha\beta_2}),\\
		\frac{c_\bullet}{\eeuler}\left( \frac{Q_{{\lambda_\alpha}}\overline{Q}_{\lambda_\alpha}}{t_{\alpha\beta_2}}\right) &= \prod_{k,k', l,l'} f((k-k')s_{\alpha\beta_1} + (l-l'-1)s_{\alpha\beta_2}),\\
		\frac{c_\bullet}{\eeuler}\left( \frac{Q_{{\lambda_\alpha}}\overline{Q}_{\lambda_\alpha}}{t_{\alpha\beta_1}}\right) &= \prod_{k,k',l,l'} f((k-k'-1)s_{\alpha\beta_1} + (l-l')s_{\alpha\beta_2}),\\
		\frac{c_\bullet}{\eeuler}\left( Q_{\lambda_\alpha}\overline{Q}_{\lambda_\alpha}\right) \;\,&= \prod_{k,k',l,l'} f((k-k')s_{\alpha\beta_1} + (l-l')s_{\alpha\beta_2}).
	\end{align}
	Hence, $C_3$ is given by
	\begin{multline}\label{eq:C3}
		\frac{c_\bullet}{\eeuler}\left( \frac{(1-t_{\alpha\beta_1})(1-t_{\alpha\beta_2})}{t_{\alpha\beta_1}t_{\alpha\beta_2}} Q_{\lambda_\alpha}\overline{Q}_{\lambda_\alpha}\right)= \prod_{k,k',l,l'} \left(\frac{f((k-k'-1)s_{\alpha\beta_1} + (l-l'-1) s_{\alpha\beta_2})}{f((k-k')s_{\alpha\beta_1} + (l-l'-1)s_{\alpha\beta_2})} \right.\\\left. \cdot \frac{f((k-k')s_{\alpha\beta_1} + (l-l')s_{\alpha\beta_2})}{f((k-k'-1)s_{\alpha\beta_1} + (l-l')s_{\alpha\beta_2})}\right).
	\end{multline}

	\paragraph{Computing $\frac{c_\bullet}{\eeuler}(E_{\alpha\beta_1})$.}
	The total Chern/Euler class of $\mathsf{E}_{\alpha\beta_1}$ is given by 
	\begin{align}
		\frac{c_\bullet}{\eeuler}(\mathsf{E}_{\alpha\beta_1}) &= \frac{c_\bullet}{\eeuler}\left(t_{\alpha\beta_1}^{-1} \frac{\mathsf{F}_{\alpha\beta_1}(t_{\alpha\beta_2})}{1-t_{\alpha\beta_1}^{-1}} - \frac{\mathsf{F}_{\alpha\beta_1}(t_{\alpha\beta_2}t_{\alpha\beta_1}^{-m_{\alpha\beta_1}})}{1-t_{\alpha\beta_1}^{-1}}\right) \\
		&= \frac{c_\bullet}{\eeuler} \Bigg(\frac{1}{1-t_{\alpha\beta_1}^{-1}}\Biggl( \underbrace{\spvertund{3.5ex}t_{\alpha\beta_1}^{-1}\mathsf{F}_{\alpha\beta_1}(t_{\alpha\beta_2}) - \mathsf{F}_{\alpha\beta_1}(t_{\alpha\beta_2}t_{\alpha\beta_1}^{-m_{\alpha\beta_1}})}_{\sum_{l=1}^{d_{\alpha\beta_1}}t_{\alpha\beta_2}^{-l} (t_{\alpha\beta_1}^{-1} - t_{\alpha\beta_1}^{l\cdot m_{\alpha\beta}})} \Biggr)\Bigg).
	\end{align}
	The term $t_{\alpha\beta_1}^{-1} - t_{\alpha\beta_1}^{l\cdot m_{\alpha\beta}}$ has a zero at $t_{\alpha\beta_1} = 1$, so we can factorise $1-t_{\alpha\beta_1}^{-1}$ out:
	Observe that
	\begin{equation}
		(1-t_{\alpha\beta_1}^{-1}) \sum_{l=1}^{d_{\alpha\beta_1}} t_{\alpha\beta_2}^{-l} \left( \sum_{k=0}^l t_{\alpha\beta_1}^k \right) = \sum_{l=1}^{d_{\alpha\beta_1}} \left( t_{\alpha\beta_2}^{-l}t_{\alpha\beta_1}^l - t_{\alpha\beta_2}^{-l} t_{\alpha\beta_1}^{-1}\right).
	\end{equation}
	Therefore, 
	\begin{align}
		\sum_{l=1}^{d_{\alpha\beta_1}} t_{\alpha\beta_2}^{-l} \left( \sum_{k=0}^l t_{\alpha\beta_1}^k \right) &=  \frac{-\sum_{l=1}^{d_{\alpha\beta_1}}t_{\alpha\beta_1}^{-1} t_{\alpha\beta_2}^{-l}}{1-t_{\alpha\beta_1}^{-1}} - \frac{-\sum_{l=1}^{d_{\alpha\beta_1}} (t_{\alpha\beta_2}t_{\alpha\beta_1}^{-1})^{-l}}{1-t_{\alpha\beta_1}^{-1}} \\
		&= t_{\alpha\beta_1}^{-1} \frac{\mathsf{F}_{\alpha\beta_1}(t_{\alpha\beta_2})}{1-t_{\alpha\beta_1}^{-1}} - \frac{\mathsf{F}_{\alpha\beta_1}(t_{\alpha\beta_2}t_{\alpha\beta_1}^{-1})}{1-t_{\alpha\beta_1}^{-1}} .
	\end{align}
	Hence,
	\begin{equation}\label{eq:C'}
		\frac{c_\bullet}{\eeuler}(\mathsf{E}_{\alpha\beta_1}) = \frac{c_\bullet}{\eeuler} \left(\sum_{l=1}^{d_{\alpha\beta_1}} t_{\alpha\beta_2}^{-l} \left( \sum_{k=0}^l t_{\alpha\beta_1}^k \right)  \right) = \prod_{l=1}^{d_{\alpha\beta_1}} \prod_{k=0}^l f(ks_{\alpha\beta_1} - l s_{\alpha\beta_2}). 
	\end{equation}
	
	\paragraph{Putting it All Together.}
	
	The equations \eqref{eq:Falpha}, \eqref{eq:C1}, \eqref{eq:C2}, \eqref{eq:C3}, and \eqref{eq:C'} imply
	\begin{multline}\label{eq:chern class}
		\frac{c_\bullet}{\eeuler}\left(T^\vir\Big|_{\vec{\lambda}\subset\vec{d}}\right) = \prod_\alpha \left(\prod_{(k,l)\in\lambda_\alpha}   \frac{f(ks_{\alpha\beta_1} + l s_{\alpha\beta_2})f(-(k+1)s_{\alpha\beta_1}- (l+1)s_{\alpha\beta_2})}{f((k-d_{\alpha\beta_1})s_{\alpha\beta_1} + (l-d_{\alpha\beta_2}) s_{\alpha\beta_2})}\right.\\ 
		\cdot  \prod_{k,k',l,l'} \frac{f((k-k'-1)s_{\alpha\beta_1} + (l-l'-1) s_{\alpha\beta_2})f((k-k')s_{\alpha\beta_1} + (l-l')s_{\alpha\beta_2})}{f((k-k')s_{\alpha\beta_1} + (l-l'-1)s_{\alpha\beta_2})f((k-k'-1)s_{\alpha\beta_1} + (l-l')s_{\alpha\beta_2})} \\
		\cdot \left.\prod_{k=1}^{d_{\alpha\beta_1}} \prod_{l=1}^{d_{\alpha\beta_2}} \frac{1}{f(-ks_{\alpha\beta_1} - l s_{\alpha\beta_2})}\right)\cdot \prod_{\alpha\beta_1} \prod_{l=1}^{d_{\alpha\beta_1}} \prod_{k=0}^l f(ks_{\alpha\beta_1} - l s_{\alpha\beta_2}). 
	\end{multline}
	\noindent
	
	\section{Computing Partition Functions}\label{sec:examples}
	
	\subsection{Global Equivariant Partition Function}
	
	Computing the global equivariant partition function, that is, the generating series of the virtual Euler characteristic can be done by reducing to the local partition function as follows:
	
	\begin{align}\label{eq:global parition function is made up of local}
		\tilde{Z}_\beta^S(q\mid s_1,s_2) &= \sum_{m\geq 0} q^m \int_{\left[S_\beta^{[0,m]}\right]^\vir} c_\bullet(T^\vir) \in\mathbb{Q} (s_1,s_2)[[q]]\\
		&= \sum_{[\vec{d}] =\beta} \sum_{\vec{\lambda}\subset \vec{d}} q^{|\vec{\lambda}|} \frac{c_\bullet}{\eeuler} \left(T^\vir\Big|_{\vec{\lambda}\subset\vec{d}}\right)\\
		&= \sum_{[\vec{d}] =\beta} \prod_{\alpha\beta} \frac{c_\bullet}{\eeuler} (\mathsf{E}_{\alpha\beta}) \prod_\alpha   \frac{c_\bullet}{\eeuler} (\mathsf{F}_\alpha) \sum_{\vec{\lambda}\subset\vec{d}} q^{|\vec{\lambda}|} \prod_\alpha \frac{c_\bullet}{\eeuler} (\mathsf{G}_\alpha)\\
		&=  \sum_{[\vec{d}] =\beta} \prod_{\alpha\beta} \frac{c_\bullet}{\eeuler} (\mathsf{E}_{\alpha\beta}) \prod_{\alpha} \frac{c_\bullet}{\eeuler}(\mathsf{F}_\alpha) \prod_\alpha Z_{(d_{\alpha\beta_1},d_{\alpha\beta_2})}^{U_\alpha} (q\mid s_{\alpha\beta_1},s_{\alpha\beta_2}),
	\end{align}
	where 
	\begin{equation}\label{eq:local partition function}
		\tilde{Z}^{\C^2}_{(d_1,d_2)} (q\mid s_{1},s_{2}) \coloneqq \sum_{\lambda\subset (d_{1}, d_{2})} q^{|\lambda|} \frac{c_\bullet}{\eeuler} \left( \mathsf{G}_{\lambda} \right).
	\end{equation}
	When we write $\lambda\subset (d_{1},d_{2})$, we mean that the partition $\lambda$ is contained in  $(d_1,d_2)$-rays as defined in Section \ref{sec:Od Ualpha}.


	We now consider the projective plane $S = \mathbb{P}^2$ and want to explicitly determine 
	\begin{equation}
		\eeuler^\vir (P_{n}(\mathbb{P}^2,d)) = \int_{[P_{n}(\mathbb{P}^2,d)]^\vir} c_\bullet\left(T^\vir\right). 
	\end{equation}
	To this end, we first consider the local partition functions \eqref{eq:local partition function}. 
	
	\subsection{Examples of Local Partition Functions and Non-Rationality}
	Let $p_\alpha$ be a fixed point of the $(\C^\times)^2$-action. Then the coordinates of $\T$ are given by $t_{1}, t_{2}$.  
	Denote the class of $t_{i}$ by $s_{i}$. 
	
	\subsubsection{When $(d_1,d_2) = (1,0)$}\label{sec:(0,1)}
	
	The Young diagram corresponding to a partition $\lambda$ restricted by $(d_1,d_2) = (1,0)$ of size $m = b$ looks as follows:
	\begin{figure}[H]
		\begin{equation}
			\YFrench
			\Yboxdim{1cm}
			\begin{tikzpicture}[scale=.7]
				\tyoung(0.5cm,0cm,~x{x^2}~~~~~,~{xy}~~~~~,~~~~,~~~,~~,~)
				\Yfillcolor{blue}
				\Yfillopacity{.1}
				\tyoung(0.5cm,0cm,1x{x^2}~~~~~~)
				\Yfillcolour{green}
				\Yfillopacity{.3}
				\tyoung(0.5cm,0cm,1x{x^2}~)
				\node [fill, draw, circle, minimum width=3pt, inner sep=0pt] at (0.5,1) {};
				\node[left=2pt of {(0.5,1)}, outer sep=2pt] {$(0,1)$};
				\node [fill, draw, circle, minimum width=3pt, inner sep=0pt] at (4.5,0) {};
				\node[below =2pt of {(4.5,-0.25)}, outer sep=2pt] {$(b,0)$};
				\draw [swap,decorate,decoration = {brace}] (4.5,-0.25) -- node[below=6pt] {$b$}  (0.5,-0.25); 
				\node[above = 2pt of {(4.5,4)}] {$\mathcal{I}_D$};
				\node[above right = 2pt of {(6.5,2.5)}] {$\mathcal{I}_Z$};
				\draw[-, very thick](0.5,1)--(7.5,1);
				\draw[-, very thick](0.5,1)--(0.5,6);
				\draw[-, very thick](4.5,0)--(4.5,4);
				\draw[-, very thick](4.5,0)--(9.5,0);
			\end{tikzpicture}
		\end{equation}
		\caption{Case $(d_1,d_2) = (1,0), \;b=4$} \label{fig:case0,1}
	\end{figure}
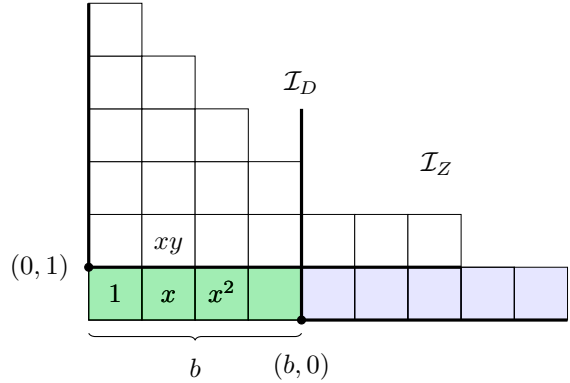
	Therefore,
	\begin{equation}
		Q_\lambda = \sum_{i=0}^{b-1} t_1^i.
	\end{equation}
	In this case the finite vertex contribution $\mathsf{G}_\alpha$ is 
	\begin{align}
		\mathsf{G}_{\lambda\subset(1,0)} &= \overline{Q}_\alpha t_{1}^{-1} t_{2}^{-1} - Q_\alpha\overline{Q}_{\alpha} (1-t_{1}^{-1})(1-t_{2}^{-1}) + Q_\alpha(1-t_{2}^{-1})\\
		&= \overline{Q}_\alpha t_{1}^{-1}t_{2}^{-1} -  Q_\alpha (t_{2}^{-1})\overline{Q}_\alpha (1-t_{1}^{-1} - 1)\\
		&= \overline{Q}_\alpha t_{1}^{-1}t_{2}^{-1} + Q_\alpha t_{1}^{-b} (1-t_{2}^{-1})\\
		&= \sum_{i=1}^b t_{1}^{-i},
	\end{align}
	where $m=b$ is the size of the partition $\lambda$. 
	This yields 
	\begin{equation}\label{eq:(0,1)}
		\tilde{Z}_{(0,1)}^{\C^2} (q\mid s_{1},s_{2}) = \sum_{b\geq 0}q^b \prod_{i=1}^b \frac{1-is_{1}}{-is_{1}}  = (1-q)^{\frac{1}{s_1} - 1} \notin \mathbb{Q}(s_{1})(q),
	\end{equation}
	which shows that the local partition function is in fact not rational.
	\begin{rem}
		Note that 
		\begin{equation}
			\tilde{Z}_{(0,1)}^{\C^2} (q\mid s_{1},s_{2}) = \tilde{Z}_{(1,0)}^{\C^2} (q\mid s_{2}, s_{1}).
		\end{equation}
	\end{rem}
	
	\subsubsection{When $(d_1,d_2) = (1,1)$}\label{sec:(1,1)}
	
	We want to compute 
	\begin{align}
		\mathsf{G}_{\lambda \subset (1,1)} &= (1-t_{1}^{-1}t_{2}^{-1})Q_\lambda + \frac{\overline{Q}_\lambda}{t_{1}t_{2}} - (1-t_{1}^{-1})(1-t_{2}^{-1}) Q_\lambda \overline{Q}_\lambda\\
		&= \underbrace{(1-t_{1}^{-1}t_{2}^{-1})Q_\lambda }_{\eqqcolon A} \underbrace{+ t_{1}^{-1} t_{2}^{-1} \overline{Q}_\lambda}_{\eqqcolon B} \underbrace{- (1-t_{1}^{-1} -t_{2}^{-1} + t_{1}^{-1} t_{2}^{-1})Q_\lambda \overline{Q}_\lambda}_{\eqqcolon C},
	\end{align}
	where the partition $\lambda$ is of size $m=a+b+1$ and its Young diagram is of the form 
	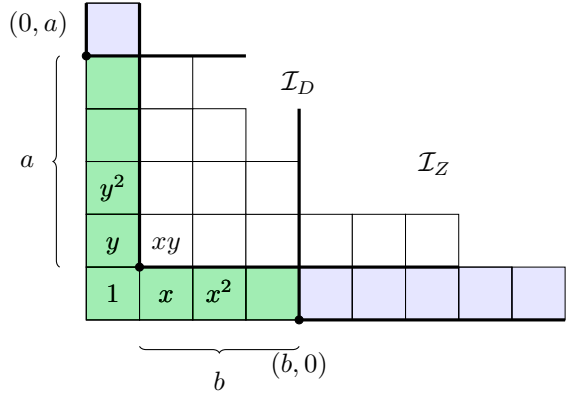
\begin{figure}[H]
		\begin{equation}
			\YFrench
			\Yboxdim{1cm}
			\begin{tikzpicture}[scale=.7]
				\tyoung(0.5cm,0cm,~x{x^2}~~~~~,~{xy}~~~~~,~~~~,~~~,~~,~)
				\Yfillcolor{blue}
				\Yfillopacity{.1}
				\tyoung(0.5cm,0cm,1x{x^2}~~~~~~,y,{y^2},~,~,~)
				\Yfillcolour{green}
				\Yfillopacity{.3}
				\tyoung(0.5cm,0cm,1x{x^2}~,y,{y^2},~,~)
				\node [fill, draw, circle, minimum width=3pt, inner sep=0pt] at (1.5,1) {};
				\node [fill, draw, circle, minimum width=3pt, inner sep=0pt] at (4.5,0) {};
				\node[below =2pt of {(4.5,-0.25)}, outer sep=2pt] {$(b,0)$};
				\node [fill, draw, circle, minimum width=3pt, inner sep=0pt] at (0.5,5) {};
				\node[above left=2pt of {(.5,5)}, outer sep=2pt] {$(0,a)$};
				\draw [decorate,decoration = {brace}] (0,1) -- node[left=6pt] {$a$}  (0,5); 
				\draw [swap,decorate,decoration = {brace}] (4.5,-0.5) -- node[below=6pt] {$b$}  (1.5,-0.5); 
				\node[above = 2pt of {(4.5,4)}] {$\mathcal{I}_D$};
				\node[above right = 2pt of {(6.5,2.5)}] {$\mathcal{I}_Z$};
				\draw[-, very thick](1.5,1)--(7.5,1);
				\draw[-, very thick](1.5,1)--(1.5,6);
				\draw[-, very thick](4.5,0)--(4.5,4);
				\draw[-, very thick](4.5,0)--(9.5,0);
				\draw[-, very thick](.5,5)--(3.5,5);
				\draw[-, very thick](.5,5)--(.5,6);
			\end{tikzpicture}
		\end{equation}
		\caption{Case $(d_1,d_2) = (1,1), \;a=5,\;b=4$} \label{fig:case1,1}
	\end{figure}
	In this case, 
	\begin{equation}
		Q_\lambda = 1 + \sum_{i=1}^{b} t_1^i + \sum_{j=1}^a t_2^j.
	\end{equation}
	Thus,
	\begin{multline}
		A+B = \sum_{i=1}^b t_{1}^i + \sum_{j=1}^a t_{2}^j + 1 - t_{2}^{-1} \sum_{i=0}^{b-1} t_{1}^i - t_{1}^{-1} \sum_{j=0}^{a-1} t_{2}^{j} + t_{2}^{-1} \sum_{i=1}^b t_{1}^{-(i+1)} + t_{1}^{-1} \sum_{j=1}^a t_{2}^{-(j+1)}.
	\end{multline}
	Furthermore,
	\begin{multline}
		Q_\lambda \overline{Q}_\lambda = t_{1} \frac{t_{1}^b-1}{t_{1}-1} \sum_{i=1}^b t_{1}^{-i} + \frac{(t_{1}^b-1)(t_{2}^a-1)}{(t_{1}-1)(t_{2}-1)} t_{1}t_{2}^{-a} + \sum_{i=1}^b t_{1}^i 
		+ \frac{(t_{1}^b-1)(t_{2}^a-1)}{(t_{1}-1)(t_{2}-1)} t_{1}^{-b}t_{2} \\+ t_{2} \frac{ t_{2}^a -1}{ t_{2}-1} \sum_{j=1}^a  t_{2}^{-j} 
		+ \sum_{j=1}^a  t_{2}^j + \sum_{i=1}^b  t_{1}^{-i} + \sum_{j=1}^a  t_{2}^{-j} +1.
              \end{multline}
              Simplifying the expression, we get
	\begin{equation}
		A+B+C = t_{1}^{-(b+1)} t_{2}^a + t_{1}^b t_{2}^{-(a+1)} - t_{1}^{-1}t_{2}^{-1} + \sum_{i=1}^b t_{1}^{-i} + \sum_{j=1}^a t_{2}^{-j}.
	\end{equation}
	Thus, the local partition function is given by 
	\begin{multline}
		\tilde{Z}_{(1,1)}^{\C^2} (q\mid s_{1}, s_{2}) = \frac{1-s_{1} - s_{2}}{-s_{1} -s_{2}} \cdot \Bigg( 1 + \sum_{a,b\geq 0} q^{a+b+1} \prod_{i=1}^b \left(\frac{1-is_{1}}{-is_{1}}\right) \prod_{j=1}^a \left(\frac{1-js_{2}}{-js_{2}}\right)\\ \frac{1+bs_{1}-(a+1)s_{2}}{bs_{1}-(a+1)s_{2}}\frac{1-(b+1)s_{1} +as_{2}}{-(b+1)s_{1} + a s_{2}}\frac{-s_{1}-s_{2}}{1-s_{1}-s_{2}} \Bigg).
	\end{multline}
	This the Laurent expansion of the function 
	\begin{equation}\label{eq:(1,1)}
		\tilde{Z}_{(1,1)}^{\C^2} (q\mid s_1,s_2) = \frac{1-s_1-s_2}{-s_1-s_2} \left( (1-q)^{\frac{1}{s_1} + \frac{1}{s_2}} \left( 1+ \frac{s_1+s_2}{-1+s_1+s_2} \frac{q}{(1-q)^2}\right) \right).
	\end{equation}

	\subsubsection{When $(d_1,d_2) = (2,0)$}\label{sec:(0,2)}
	
	For $(d_{1},d_{2}) = (2,0)$ and a partition $\lambda$ of size $2a+b$, we have 
	\begin{equation}
		Q_\lambda = (1+t_2)\sum_{i=0}^{a-1} t_{1}^i + \sum_{j=0}^{b-1} t_{1}^{a+j}.
	\end{equation}
	The corresponding Young diagram is:
	\begin{figure}[H]
		\begin{equation}
			\YFrench
			\Yboxdim{1cm}
			\begin{tikzpicture}[scale=.7]
				\tyoung(0.5cm,0cm,~x{x^2}~~~~~,~{xy}~~~~~,~~~~,~~~,~~,~)
				\Yfillcolor{blue}
				\Yfillopacity{.1}
				\tyoung(0.5cm,0cm,1x{x^2}~~~~~~,y{xy}~~~~~)
				\Yfillcolour{green}
				\Yfillopacity{.3}
				\tyoung(0.5cm,0cm,1x{x^2}~~~,y{xy}~~)
				\node [fill, draw, circle, minimum width=3pt, inner sep=0pt] at (0.5,2) {};
				\node[left=2pt of {(0.5,2)}, outer sep=2pt] {$(0,2)$};
				\node [fill, draw, circle, minimum width=3pt, inner sep=0pt] at (4.5,1) {};
				\node [fill, draw, circle, minimum width=3pt, inner sep=0pt] at (6.5,0) {};
				\draw [swap,decorate,decoration = {brace}] (4.5,-0.5) -- node[below=9.5pt] {$a$}  (0.5,-0.5); 
				\draw [swap,decorate,decoration = {brace}] (6.5,-0.5) -- node[below=6pt] {$b$}  (4.5,-0.5); 
				\node[above = 2pt of {(4.5,4)}] {$\mathcal{I}_D$};
				\node[above right = 2pt of {(6.5,2.5)}] {$\mathcal{I}_Z$};
				\draw[-, very thick](0.5,2)--(7.5,2);
				\draw[-, very thick](0.5,2)--(0.5,6);
				\draw[-, very thick](4.5,1)--(4.5,4);
				\draw[-, very thick](4.5,1)--(8.5,1);
				\draw[-, very thick](6.5,0)--(6.5,3);
				\draw[-, very thick](6.5,0)--(9.5,0);
			\end{tikzpicture}
		\end{equation}
		\caption{Case $(d_1,d_2) = (2,0), \;a=4, \;b=2$} \label{fig:case0,2}
	\end{figure}
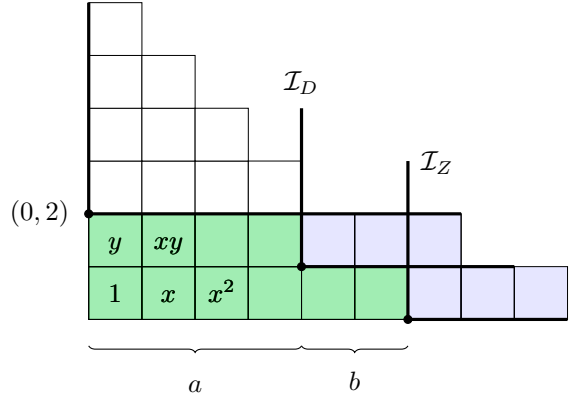
	We want to compute
	\begin{equation}\label{eq:Galpha for d1,d2 = 0,2}
		\mathsf{G}_{\lambda \subset (2,0)} = \underbrace{Q_\lambda (1-t_{2}^{-2})}_{\eqqcolon	A} \underbrace{+ \frac{\overline{Q}_\lambda}{t_{1}t_{2}}}_{\eqqcolon B} \underbrace{- (t_{1}^{-1}-1)(t_{2}^{-1}-1) Q_\lambda \overline{Q}_\lambda}_{\eqqcolon C}.
	\end{equation}
	Observe that  
	\begin{align}\label{eq:A}
		A &= (1+t_{2})(1-t_{2}^{-2}) \sum_{i=0}^{a-1}t_{1}^i + (1-t_{2}^{-2}) \sum_{j=0}^{b-1} t_{1}^{a+j}\\
		&= (1-t_{2}^{-1})(1+t_{2}^{-1})(1+t_{2})\sum_{i=0}^{a-1} t_{1}^i + (1-t_{2}^{-2})\sum_{j=0}^{b-1}t_{1}^{a+j}
	\end{align}
	and
	\begin{align}\label{eq:B}
		B &= (t_{2}^{-1}+t_{2}^{-2})\sum_{i=1}^a t_{1}^{-i} + \sum_{j=1}^b t_{1}^{-(a+j)}t_{2}^{-1}\\
		&= t_{1}^{-1}\sum_{j=1}^{b} t_{1}^{-j} + t_{2}^{-1}\sum_{i=1}^a t_{1}^{-(i+b)} + t_{2}^{-2}\sum_{i=1}^a t_{1}^{-i}. 
	\end{align}
	We are left to determine $C$. Note that 
	\begin{align}\label{eq:computing (0,2): VVbar} 
		Q_\lambda \overline{Q}_\lambda &= \left((1+t_{2}) \sum_{i=0}^{a-1}t_{1}^i + \sum_{j=0}^{b-1}t_{1}^{a+j}\right)\left((1+t_{2}^{-1}) \sum_{i=0}^{a-1}t_{1}^{-i}+ \sum_{j=0}^{b-1}t_{1}^{-(a+j)}\right)\\
		&= (1+t_2)(1+t_2^{-1}) \sum_{i=0}^{a-1} t_1^i \sum_{i'=0}^{a-1}t_1^{-i'} + (1+t_2)\sum_{i=0}^{a-1}t_1^i \sum_{j=0}^{b-1} t_1^{-(a+j)} \\ & \qquad \qquad \qquad + (1+t_2^{-1})\sum_{j=0}^{b-1} t_1^{a+j} \sum_{i=0}^{a-1} t_1^{-i} + \sum_{j=0}^{b-1} t_1^{a+j} \sum_{j'=0}^{b-1} t_1^{-(a+j')}\\
		&= (1+t_2)(1+t_2^{-1}) \sum_{i,i'=0}^{a-1} t_1^{i-i'} + (1+t_2)\sum_{i=0}^{a-1}\sum_{j=0}^{b-1}t_1^{i-(a+j)} \\ & \qquad \qquad\qquad \qquad \qquad + (1+t_2^{-1})\sum_{j=0}^{b-1}  \sum_{i=0}^{a-1} t_1^{a+j-i} + \sum_{j,j'=0}^{b-1} t_1^{j-j'}.
	\end{align}
	Denote 
	\begin{equation}
		C' \coloneqq - (t_{1}^{-1}-1)(t_{2}^{-1} - 1)(1+t_{2})(1+t_{2}^{-1}) \sum_{i,i' = 0}^{a-1} t_{1}^{i-i'}  - (t_{1}^{-1}-1)(t_{2}^{-1} - 1) \sum_{j,j'=0}^{b-1} t_{1}^{j-j'}
	\end{equation}
	and 
	\begin{multline}
		C'' \coloneqq - (t_{1}^{-1}-1)(t_{2}^{-1} - 1)(1+t_{2})\sum_{i=0}^{a-1}\sum_{j=0}^{b-1}t_{1}^it_{1}^{-(a+j)}
		- (t_{1}^{-1}-1)(t_{2}^{-1} - 1) (1+t_{2}^{-1})\sum_{i=0}^{a-1} \sum_{j=0}^{b-1} t_{1}^{-i}t_{1}^{a+j},
	\end{multline}
	so that $C = C' + C''$. 
	Observe that 
	\begin{equation}
		\sum_{i,i' = 0}^{a-1} t_{1}^{i-i'} =  \frac{t_{1}^{-(a-1)}(t_{1}^a-1)}{t_{1}-1} \cdot \sum_{i=0}^{a-1}t_{1}^i = \frac{t_{1} - t_{1}^{-(a-1)}}{t_{1}-1} \cdot \sum_{i=0}^{a-1}t_{1}^i .
	\end{equation}
	There is an analogous formula for the sum over $j$ and $j'$. And hence,
	\begin{align}\label{eq:computing (0,2): VVbar single terms}
		C' &=-(1-t_{2}^{-1})(1+t_{2})(1+t_{2}^{-1}) (1-t_{1}^{-a}) \sum_{i=0}^{a-1} t_{1}^i - (1-t_{2}^{-1})(1 - t_{1}^{-b}) \sum_{j=0}^{b-1}t_{1}^j\\
		&= -(1-t_{2}^{-2})(1+t_{2})\left(\sum_{i=0}^{a-1}t_{1}^i - \sum_{i=1}^a t_{1}^{-i}\right)- (1-t_{2}^{-1}) \left(\sum_{j=0}^{b-1} t_{1}^j - \sum_{j=1}^b t_{1}^{-j}\right) .
	\end{align}
	For $C''$ we have,
	\begin{align}\label{eq:computing (0,2): VVbar mixed terms}
		C'' &= -(t_{2}-t_{2}^{-1})(t_{1}^a-1)t_{1}^{-(a+b)}\sum_{j=0}^{b-1} t_{1}^j - (1-t_{2}^{-2})(t_{1}^a - 1)\sum_{j=0}^{b-1}t_{1}^j.
	\end{align}
	Combining \eqref{eq:computing (0,2): VVbar single terms} and \eqref{eq:computing (0,2): VVbar mixed terms}, and the expressions \eqref{eq:A} and \eqref{eq:B} we have for $A$ and $B$, we get  
	\begin{align}
		\mathsf{G}_{\lambda \subset (2,0)} &= (1+t_2) \sum_{i=1}^a t_1^{-i} + (t_2^{-1}-t_2^{-2}) \sum_{j=0}^{b-1} t_1^j + (1-t_2) \sum_{j=1}^b t_1^{-j} + t_2 \sum_{j=1}^b t_1^{-(a+j)}\\
        &= \sum_{i=1}^a t_1^{-i} + (t_2^{-1}-t_2^{-2}) \sum_{j=0}^{b-1} t_1^j + \sum_{j=1}^b t_1^{-j}  + t_2 \sum_{i=1}^a t_1^{-(b+i)}.
	\end{align}
	Hence, 
	\begin{multline}\label{eq:(0,2)}
		\tilde{Z}_{(2,0)}^{\C^2} (q\mid s_{1},s_{2}) = \sum_{a,b\geq 0} q^{2a+b} \left( \prod_{i=1}^a \frac{1-is_{1}}{-is_{1}}\frac{1-(i+b)s_{1}+s_{2}}{-(i+b)s_{1}+s_{2}} \right.\\\left.\prod_{j=1}^b \frac{1-js_{1}}{-js_{1}} \prod_{i=0}^{b-1} \frac{1+is_{1} -s_{2}}{is_{1}-s_{2}}\frac{is_{1}-2s_{2}}{1+is_{1}-2s_{2}} \right).
	\end{multline}
	We do not have a closed formula for a function with this Laurent expansion, but we suspect that it exists. 
	\subsection{Python Program}\label{sec:program}
	The website of the author contains the source code of a program, which computes 
	\begin{equation}
		\int_{[P_{n}(\mathbb{P}^2,d)]^\vir} c_\bullet(T^\vir)
	\end{equation}
	for any given $n$ and $d$ using the method described above. Figures \ref{fig:numbers1} and \ref{fig:numbers2} contain values for $d=1,2,3,4$ and low $m$, where as in Section \ref{sec:nested hilbert schemes}, 
	\begin{equation}
		m = n+ \frac{d(d-3)}{2}.
	\end{equation}
 
	\begin{figure}[H]
		\begin{center} 
			\begin{tabular}{l|l|l}
				$d$ & $m$ & output$\qquad\quad\quad\;\;\,\,$\\
				\hline 
				1 & 0 &  3\\
				1 & 1 & 6\\
				1 & 2 & 9\\
				1 & 3 & 12\\
				1 & 4 & 15\\
				1 & 5 & 18\\
				1 & 6 & 21\\
				1 & 7 & 24\\
				1 & 8 & 27\\
				1 & 9 & 30\\
				1 & 10 & 33\\
				1 & 11 & 36\\
				1 & 12 & 39\\
				1 & 13 & 42\\
				1 & 14 & 45\\
				1 & 15 & 48\\
				1 & 16 & 51\\
				1 & 17 & 54\\
				1 & 18 & 57\\
				1 & 19 & 60\\
				1 & 20 & 63\\
				1 & 21 & 66\\
				1 & 22 & 69\\
				1 & 23 & 72\\
				1 & 24 & 75
			\end{tabular}
			\begin{tabular}{l|l|l}
				$d$ & $m$ & output$\qquad\qquad\quad\quad\,$\\
				\hline 
				2 & 0 & 6\\
				2 & 1 & 15\\
				2 & 2 & 36\\
				2 & 3 & 66\\
				2 & 4 & -336\\
				2 & 5 & -8019\\
				2 & 6 &-70098\\
				2 & 7 & -399804\\
				2 & 8 & -1740870\\
				2 & 9 & -6260277\\
				2 & 10 & -19487496\\
				2 & 11 &  -54159930\\
				2 & 12 & -137321340\\
				2 & 13 & -322688919\\
				2 & 14 & -711195678\\
				2 & 15 & -1483772280\\
				2 & 16 & -2951904786\\
				2 & 17 & -5633344377\\
				2 & 18 & -10362608436\\
				2 & 19 & -18448659894\\
				2 & 20 & -31895447976\\
				2 & 21 & -53704906971\\
				2 & 22 & -88286612970\\
				2 & 23 & -142003668276\\
				2 & 24 & -223890600030
			\end{tabular}
            \caption{Virtual Euler characteristic for $d=1,2$ and low $m$}
		      \label{fig:numbers1}
		\end{center}
    \end{figure}
    \begin{figure}[H]
		\begin{center}
			\begin{tabular}{l|l|l}
				$d$ & $m$ & output$\qquad\quad\;\,$\\
				\hline 
				3 & 0 & 10\\
				3 & 1 & 27\\
				3 & 2 & 72\\
				3 & 3 & 154\\
				3 & 4 & 306\\
				3 & 5 & -19737\\
				3 & 6 & -1349404\\
				3 & 7 & -32053869\\
				3 & 8 & -430135668\\
				3 & 9 & -3946790877\\
				3 & 10 & -27473408784\\
				3 & 11 & -154768875579\\
				3 & 12 & -736999029842\\
				3&  13 & -3059890203483\\
				3&  14 & -1133301372836 \\
				3&  15 & -38104196925509\\
				3&  16 & -117902025110844\\
			\end{tabular}
			\begin{tabular}{l|l|l}
				$d$ & $m$ & output$\qquad\qquad\quad\quad\,$\\
				\hline 
				4 & 0 & 15\\
				4 & 1 & 42\\
				4 & 2 & 117\\
				4 & 3 & 264\\
				4 & 4 & 561\\
				4 & 5 & 1080\\
				4 & 6 & 26058\\
				4 & 7 & 16548006\\
				4 & 8 & 1842925419\\
				4 & 9 & 80399046090\\
				4 & 10 & 1942340199207\\
				4 & 11 & 30960585072144\\
				4 & 12 & 361026356454855\\
				4 & 13 & 3293495920441878\\
				4 & 14 & 24626906563808097\\
				4 & 15 & 156153491429509728\\
				4 & 16 & 861447562288733412
			\end{tabular}
		\end{center}
		\caption{Virtual Euler characteristic for $d=3,4$ and low $m$}
		\label{fig:numbers2}
	\end{figure}
	
	Note that when $d=1$ and for the cases where $d$ is general and $n\leq d+1$, we can determine the Euler characteristic via different methods. 
	\paragraph{When $d=1$.}
	For $S=\mathbb{P}^2$, consider the forgetful map 
	\begin{equation}
		S^{[0,m]}_d \to \operatorname{Hilb}_d, (Z,D) \mapsto D.
	\end{equation}
	As $d=1$, we have $D\cong\mathbb{P}^1$. Thus, the Hilbert scheme of divisors with $d=1$ parametrises lines in $\mathbb{P}^2$. Hence, $ \operatorname{Hilb}_{d=1} = \left(\mathbb{P}^2\right)^\lor$ and the fibres of this map are given by the Hilbert scheme of points 
	\begin{equation}
		D^{[m]} \cong \mathbb{P}^m.
	\end{equation}
	Thus, the Euler characteristic of $S^{[0,m]}_{d=1}$ is given by 
	\begin{equation}
		\chi\left(\left(\mathbb{P}^2\right)^\lor\right) \cdot \chi(\mathbb{P}^m) = 3 \cdot (m+1).
	\end{equation}
	\paragraph{When $n\leq d+1$, $d\geq1$.}
	Consider the other forgetful map 
	\begin{equation}
		S^{[0,m]}_{d} \to S^{[m]}, (Z,D)\mapsto Z.
	\end{equation}
	The fibres of this map are given by 
	\begin{equation}
		\mathbb{P}H^0(\mathcal{I}_Z(d)).
	\end{equation} 
	The higher cohomology of $\mathcal{I}_Z(d)$ vanishes. In degree $2$ this follows from the fact that 
	\begin{equation}
		\operatorname{Hom}(\mathcal{I}_Z,\mathcal{O}(-3-d))^\lor \hookrightarrow \operatorname{Hom}(\mathcal{O},\mathcal{O}(-3-d))^\lor = 0
	\end{equation}
	as $p_g =0$. In degree $1$ the vanishing follows from $n\leq d+1$. Indeed, consider the short exact sequence 
	\begin{equation}
		0\to \mathcal{I}_Z(d) \to \mathcal{O}(d) \to \mathcal{O}_Z(d) \to 0.
	\end{equation}
	This induces the long exact sequence 
	\begin{equation}
		0 \to H^0(\mathcal{I}_Z(d)) \to H^0(\mathcal{O}(d))\xrightarrow{f} H^0(\mathcal{O}_Z(d)) \to H^1(\mathcal{I}_Z(d)) \to H^1(\mathcal{O}(d)) \to 0.
	\end{equation}
	The map $f$ is surjective because $\mathcal{O}(d)$ is $d$-very ample and $n\leq d+1$. The Euler characteristic is computed as
	\begin{align}
		\chi(S^{[m]}) \chi(\mathcal{I}_Z(d)) &=   \chi(S^{[m]})\cdot \left(1 + \frac{d(d+3)}{2} - m\right),
	\end{align}
	where $\chi(S^{[m]})$ is the Euler characteristic of the Hilbert scheme of $m$ points of $\mathbb{P}^2$. This can be comptuted via G\"ottsche's formula \eqref{eq:goettsche} as the $m$th coefficient in the series 
	\begin{multline}
		1 + 3 q + 9 q^2 + 22 q^3 + 51 q^4 + 108 q^5 + 221 q^6 + 429 q^7 + 810 q^8 + 1479 q^9 + 2640 q^{10} \\
		+ 4599 q^{11} + 7868 q^{12} + 13209 q^{13} + O(q^{14}).
	\end{multline}

	\subsection{Examples of Global Partition Functions}
	
	Denote the fixed points of the $(\C^\times)^2$-action on $\mathbb{P}^2$ by $p_\alpha$, $p_\beta$, and $p_\gamma$.
	Denote by $s_1$ the class of $t_1 = t_{\alpha\beta}$ and by $s_2$ the class of $t_2 = t_{\alpha\gamma}$, the other $t$'s are noted in the diagram:
	
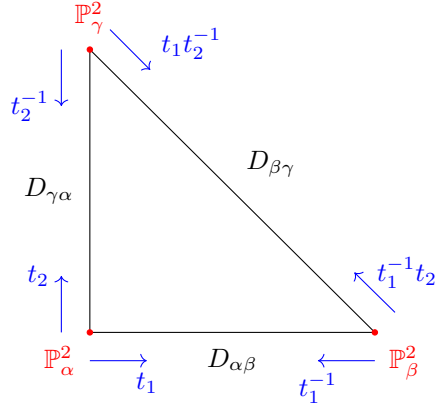
\begin{equation}
		\begin{tikzpicture}[scale=1.5]
			\draw[-] (0,0)--(2.5,0);
			\draw[-] (0,0)--(0,2.5);
			\draw[-] (2.5,0)--(0,2.5);
			\node[below left = 2.5pt of {(0,0)}] {\textcolor{red}{$\mathbb{P}^2_\alpha$}};
			\node [fill, draw, circle, minimum width=2pt, inner sep=0pt,color=red] at (0,0) {};
			\node[below right= 2.5pt of {(2.5,0)}] {\textcolor{red}{$\mathbb{P}^2_\beta$}};
			\node [fill, draw, circle, minimum width=2pt, inner sep=0pt,color=red] at (2.5,0) {};
			\node[above = 2.5pt of {(0,2.5)}] {\textcolor{red}{$\mathbb{P}^2_\gamma$}};
			\node [fill, draw, circle, minimum width=2pt, inner sep=0pt,color=red] at (0,2.5) {};
			\node[above right = 2.5pt of {(1.25,1.25)}] {\textcolor{black}{$D_{\beta\gamma}$}};
			\node[left = 2.5pt of {(0,1.25)}] {\textcolor{black}{$D_{\gamma\alpha}$}};
			\node[below = 2.5pt of {(1.25,0)}] {\textcolor{black}{$D_{\alpha\beta}$}};
			\draw[->, color=blue] (0,-.25)--(.5,-.25);
			\node[below = 1pt of {(.5,-.25)}] {\textcolor{blue}{$t_1$}};
			\draw[->, color=blue] (-.25,0)--(-.25,.5);
			\node[left = 1pt of {(-.25,.5)}] {\textcolor{blue}{$t_2$}};
			\draw[->, color=blue] (2.5,-.25)--(2,-.25);
			\node[below = 1pt of {(2,-.25)}] {\textcolor{blue}{$t_1^{-1}$}};
			\draw[->, color=blue] (2.6768,.1768)--(2.32,.53);
			\node[right = 5pt of {(2.32,.53)}] {\textcolor{blue}{$t_1^{-1}t_2$}};
			\draw[->, color=blue] (.1768,2.6768)--(.53,2.32);
			\node[above right = 1pt of {(.53,2.32)}] {\textcolor{blue}{$t_1t_2^{-1}$}};
			\draw[->, color=blue] (-.25,2.5)--(-.25,2);
			\node[left = 1pt of {(-.25,2)}] {\textcolor{blue}{$t_2^{-1}$}};
		\end{tikzpicture}
	\end{equation}
	
	\subsubsection{When $d=1$}\label{sec:d=1}
	
	The aim of this section is to prove Theorem \ref{thm:P2,d=1}, i.e., to determine the global partition function 
	\begin{equation}
		Z_{d=1}^{\mathbb{P}^2} (q) = Z_{d=1}^{\mathbb{P}^2} (q\mid s_1,s_2) \Big|_{s_1=s_2=0} =  \sum_{|\vec{d}|=1} \sum_{\vec{\lambda}\subset\vec{d}} q^{|\vec{\lambda}|} \frac{c_\bullet}{\eeuler} \left(T^\vir\Big|_{\vec{\lambda}\subset\vec{d}}\right).
	\end{equation}
	
	To compute the partition function, we use the result from Section \ref{sec:(0,1)}: Equation \eqref{eq:(0,1)} states that for any fixed point $\mathbb{P}^2_\alpha$, we have 
	\begin{equation}
		\tilde{Z}_{(1,0)}^{U_\alpha}(q\mid s_{\alpha\beta_1},s_{\alpha\beta_2}) = (1-q)^{\frac{1}{s_{\alpha\beta_1}}-1}.
	\end{equation}
	Similarly, we have 
	\begin{equation}
		\tilde{Z}_{(0,1)}^{U_\alpha}(q\mid s_{\alpha\beta_1},s_{\alpha\beta_2}) = (1-q)^{\frac{1}{s_{\alpha\beta_2}}-1}.
	\end{equation}
	Although these individually may not be rational functions, the rational parts of the exponents cancel after accounting for the contributions of the other vertices. Hence, the global partition function is rational nonetheless. 
	The edge terms are given by 
	\begin{equation}
		\mathsf{E}_{\alpha\beta_1}^{d_{\alpha\beta_1} =0}(t_{\alpha\beta_1},t_{\alpha\beta_2}) = 0
	\end{equation}
	and 
	\begin{equation}
		\mathsf{E}_{\alpha\beta_1}^{d_{\alpha\beta_1} =1}(t_{\alpha\beta_1},t_{\alpha\beta_2}) = t_{\alpha\beta_2}^{-1} + t_{\alpha\beta_2}t_{\alpha\beta_1}^{-1}
	\end{equation}
	and the infinite vertex terms $\mathsf{F}_\alpha$ are empty sums, so equal to $0$. 
	
	Therefore, as
	\begin{equation}
		\sum_{\lambda\subset (1,0)} q^{|\lambda|} \frac{c_\bullet}{\eeuler} (\mathsf{G}_\alpha) = \tilde{Z}_{(1,0)}^{\C^2}(q\mid s_{\alpha\beta},s_{\alpha\gamma}),
	\end{equation}
	we have
	\begin{align}\label{eq:global 1}
		\sum_{|\vec{d}|=1} \prod_{\alpha\beta} \frac{c_\bullet}{\eeuler} (\mathsf{E}_{\alpha\beta})& \sum_{\vec{\lambda}\subset \vec{d}} \prod_{\alpha} q^{|\vec{\lambda}|} \frac{c_\bullet}{\eeuler}(\mathsf{G}_\alpha) \\
		= &\frac{c_\bullet}{\eeuler} (t_{\gamma\beta}^{-1} t_{\gamma\alpha}+t_{\gamma\beta}^{-1}) \cdot \left(Z_{(1,0)}^{\C^2}(q\mid s_{\gamma\alpha},s_{\gamma\beta})\cdot Z_{(1,0)}^{\C^2}(q\mid s_{\alpha\gamma}, s_{\alpha\beta})\right) \\
		&+  \frac{c_\bullet}{\eeuler} (t_{\beta\alpha}^{-1} t_{\beta\gamma}+t_{\beta\alpha}^{-1}) \cdot \left(Z_{(1,0)}^{\C^2}(q\mid s_{\beta\gamma},s_{\beta\alpha})\cdot Z_{(1,0)}^{\C^2}(q\mid s_{\gamma\beta}, s_{\gamma\alpha})\right)   \\
		&+  \frac{c_\bullet}{\eeuler} (t_{\alpha\gamma}^{-1} t_{\alpha\beta}+t_{\alpha\gamma}^{-1}) \cdot \left(Z_{(1,0)}^{\C^2}(q\mid s_{\alpha\beta},s_{\alpha\gamma})\cdot Z_{(1,0)}^{\C^2}(q\mid s_{\beta\gamma}, s_{\beta\alpha})\right).
	\end{align}
	Recall the diagram
	
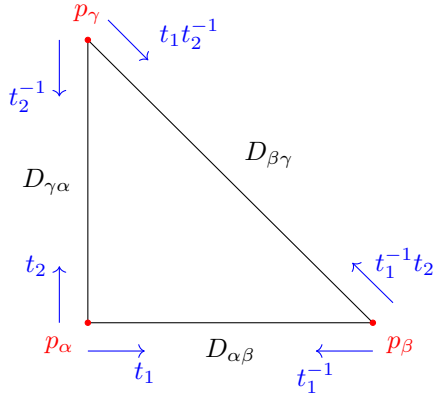
\begin{equation}
		\begin{tikzpicture}[scale=1.5]
			\draw[-] (0,0)--(2.5,0);
			\draw[-] (0,0)--(0,2.5);
			\draw[-] (2.5,0)--(0,2.5);
			\node[below left = 2.5pt of {(0,0)}] {\textcolor{red}{$p_\alpha$}};
			\node [fill, draw, circle, minimum width=2pt, inner sep=0pt,color=red] at (0,0) {};
			\node[below right= 2.5pt of {(2.5,0)}] {\textcolor{red}{$p_\beta$}};
			\node [fill, draw, circle, minimum width=2pt, inner sep=0pt,color=red] at (2.5,0) {};
			\node[above = 2.5pt of {(0,2.5)}] {\textcolor{red}{$p_\gamma$}};
			\node [fill, draw, circle, minimum width=2pt, inner sep=0pt,color=red] at (0,2.5) {};
			\node[above right = 2.5pt of {(1.25,1.25)}] {\textcolor{black}{$D_{\beta\gamma}$}};
			\node[left = 2.5pt of {(0,1.25)}] {\textcolor{black}{$D_{\gamma\alpha}$}};
			\node[below = 2.5pt of {(1.25,0)}] {\textcolor{black}{$D_{\alpha\beta}$}};
			\draw[->, color=blue] (0,-.25)--(.5,-.25);
			\node[below = 1pt of {(.5,-.25)}] {\textcolor{blue}{$t_1$}};
			\draw[->, color=blue] (-.25,0)--(-.25,.5);
			\node[left = 1pt of {(-.25,.5)}] {\textcolor{blue}{$t_2$}};
			\draw[->, color=blue] (2.5,-.25)--(2,-.25);
			\node[below = 1pt of {(2,-.25)}] {\textcolor{blue}{$t_1^{-1}$}};
			\draw[->, color=blue] (2.6768,.1768)--(2.32,.53);
			\node[right = 5pt of {(2.32,.53)}] {\textcolor{blue}{$t_1^{-1}t_2$}};
			\draw[->, color=blue] (.1768,2.6768)--(.53,2.32);
			\node[above right = 1pt of {(.53,2.32)}] {\textcolor{blue}{$t_1t_2^{-1}$}};
			\draw[->, color=blue] (-.25,2.5)--(-.25,2);
			\node[left = 1pt of {(-.25,2)}] {\textcolor{blue}{$t_2^{-1}$}};
		\end{tikzpicture}
	\end{equation}
	and that $s_i$ is the class of $t_i$. Then \eqref{eq:global 1} is equal to
	\begin{multline}
		\frac{1-s_1}{-s_1} \frac{1-s_1+s_2}{-s_1+s_2} (1-q)^{- \frac{1}{s_2} + \frac{1}{s_2} - 2}
		+ \frac{1+s_1}{s_1} \frac{1+s_2}{s_2} (1-q)^{\frac{1}{s_2-s_1} + \frac{1}{s_1-s_2} - 2}
		\\+  \frac{1-s_2}{-s_2} \frac{1+s_1 - s_2}{s_1-s_2} (1-q)^{\frac{1}{s_1} - \frac{1}{s_1} - 2},
	\end{multline}
	which simplifies to 
	\begin{equation} \label{eq:S=P2, d=1}
		Z_{d=1}^{\mathbb{P}^2}(q) = Z_{d=1}^{\mathbb{P}^2}(q\mid s_1,s_2) = \frac{3}{(1-q)^2}.
	\end{equation}
	
	\subsubsection{When $d=2$}\label{sec:d=2}
	
	We want to compute 
	\begin{equation}
		Z_{d=2}^{\mathbb{P}^2} (q) = Z_{d=2}^{\mathbb{P}^2} (q\mid s_1,s_2) \Big|_{s_1=s_2=0} =  \sum_{|\vec{d}|=2} \sum_{\vec{\lambda}\subset\vec{d}} q^{|\vec{\lambda}|} \frac{c_\bullet}{\eeuler} \left(T^\vir\Big|_{\vec{\lambda}\subset\vec{d}}\right).
	\end{equation}
	If $d=2$, then either $\vec{d}$ contains two $1$'s and one $0$ or $\vec{d}$ contains one $2$ and two $0$'s. 
	In the first case, the local partition functions are given by \eqref{eq:(0,1)} and \eqref{eq:(1,1)}, whereas in the second case, the local partition functions are given by the constant factor $1$ (if locally we have $(d_1,d_2) = (0,0)$) and \eqref{eq:(0,2)}.
	
	First, we consider the cases $\vec{d} \in \{(0,1,1), (1,0,1), (1,1,0)\}$. By \eqref{eq:(0,1)} and \eqref{eq:(1,1)}, we have 
	\begin{multline}\label{eq:(0,1,1)}
		\sum_{\vec{\lambda}\subset (0,1,1)} q^{|\vec{\lambda}|}\frac{c_\bullet}{\eeuler} (\mathsf{G}_\alpha) \frac{c_\bullet}{\eeuler} (\mathsf{G}_\beta)  \frac{c_\bullet}{\eeuler} (\mathsf{G}_\gamma) \\
		= \tilde{Z}_{(0,1)}^{\C^2} (q\mid -s_1,s_2-s_1)\cdot \tilde{Z}_{(1,1)}^{\C^2} (q\mid s_1-s_2,-s_2) \cdot \tilde{Z}_{(1,0)}^{\C^2} (q\mid s_2,s_1).
	\end{multline}
	This product is a rational function, because all the rational exponents of $(1-q)$ in \eqref{eq:(0,1)} and \eqref{eq:(1,1)} cancel out. 
	
	Second, we consider the cases $\vec{d} \in \{(0,0,2),(0,2,0),(2,0,0)\}$. By \eqref{eq:(0,2)}, we have 
	\begin{equation}\label{eq:(0,0,2)}
		\sum_{\vec{\lambda}\subset (0,0,2)} q^{|\vec{\lambda}|}\frac{c_\bullet}{\eeuler} (\mathsf{G}_\alpha) \frac{c_\bullet}{\eeuler} (\mathsf{G}_\beta)  \frac{c_\bullet}{\eeuler} (\mathsf{G}_\gamma)	= \tilde{Z}_{(2,0)}^{\C^2} (q\mid -s_2,s_1-s_2) \cdot \tilde{Z}_{(0,2)}^{\C^2} (q\mid s_1,s_2).
	\end{equation}
	If we sum all terms of the form \eqref{eq:(0,1,1)} and \eqref{eq:(0,0,2)}, we get 
	\begin{equation}
		6 + 15q + 36q^2 + 66q^3 - 336q^4 -8019q^5 -70098q^6 - \ldots
	\end{equation}
	Interpolating the first 12 coefficients of the Laurent expansion yields a conjectural formula for the coefficient of the $q^n$-term in the degree $2$ case:
	\begin{equation}
		\frac{1872360 n - 360008 n^3 + 171325 n^5 - 21714 n^7 + 1265 n^9 - 28 n^{11}}{277200}. 
	\end{equation}
	This agrees with all $25$ values in Figure \ref{fig:numbers1}. Therefore, conjecturally, the generating function is
	\begin{equation}
		\frac{6 q - 57 q^2 + 252 q^3 - 696 q^4 + 918 q^5 - 4878 q^6 + 918 q^7 - 696 q^8 + 252 q^9 - 57 q^{10} + 6 q^{11}}{(1-q)^{12}}.
	\end{equation}

	A useful criterion to determine whether a series is the expansion of a rational function is the Kronecker criterion: 
	\begin{lemma}[Kronecker]
		The power series 
		\begin{equation}
			\sum_{n\geq0} c_nz^n
		\end{equation}
		represents a rational function if and only if there exists an $m_0$ such that for all $m\geq m_0$
		\begin{equation}
			\det \begin{pmatrix}
				c_0 & c_1 & \cdots & c_m \\
				c_1 & c_2 & \cdots & c_{m+1}\\
				\vdots & \vdots & & \vdots\\
				c_m & c_{m+1} & \cdots & c_{2m}
			\end{pmatrix} =0.
		\end{equation}
	\end{lemma}
	The Kronecker criterion further supports the conjecture that $Z_{d=2}^{\mathbb{P}^2}(q)$ is the Laurent expansion of a rational function in $q$:
	Using the first 25 coefficients from Section \ref{sec:program}, we define the matrix $M_2 = (m_{ij})_{1\leq i,j\leq 13}$ by 
	\begin{equation}
		m_{ij} \coloneqq (i+j-1)\text{st} \text{ coefficient}.
	\end{equation}
	This matrix has determinant 
	\begin{equation}
		\det M_2 = 0.
	\end{equation}
	
	\subsubsection{Higher $d$}\label{sec:higher d}
	
	Based on the results in Section \ref{sec:program}, we can conjecturally deduce that 
	\begin{multline}
		Z_{d=3}^S (q) = (10 - 153 q + 1116 q^2 - 5171 q^3 + 17118 q^4 - 63495 q^5 - 898360 q^6 \\
		- 10996722 q^7 - 42987618 q^8 - 69231380 q^9 - 42987618 q^{10}\\
		- 10996722 q^{11} 	- 898360 q^{12} - 63495 q^{13} + 17118 q^{14} - 5171 q^{15} \\+ 1116 q^{16} - 153 q^{17} + 10 q^{18})(1-q)^{-18}
	\end{multline}
	and
	\begin{multline}
		Z_{d=4}^S (q) = (15 - 318 q + 3249 q^2 - 21312 q^3 + 100899 q^4 - 367596 q^5 +	1097652 q^6 \\
		+ 13378998 q^7 + 1457762922 q^8 + 40677292032 q^9 +	488193023214 q^{10} \\
		+ 2979244822788 q^{11} + 10213475593308 q^{12} + 20895483276906 q^{13} \\
		+ 26422064163513 q^{14} + 20895483276906 q^{15} + 10213475593308 q^{16}\\
		+...- 318 q^{27}+15 q^{28})(1-q)^{-24}.
	\end{multline}
	In particular, the denominators in these examples are of the form $(1-q)^{6d}$ and the numerator is of degree $d(d+3)$. Furthermore, note that the numerator is palindromic and therefore, there is a 
	\begin{equation}
		q\leftrightarrow \frac{1}{q}
	\end{equation}
	symmetry resembling the symmetry for Calabi-Yau three-folds \cite[(11)]{panda_3}.

	\printbibliography
	
\end{document}